\documentclass[12pt]{article}
\usepackage{amsmath}
\usepackage{amssymb}

\newcommand{\B}{{\mathcal B}}
\newcommand{\Ibvn}{{\mathcal IBV}^n}
\newcommand{\Ibv}{{\mathcal IBV}}
\newcommand{\qed}{\mbox{$\quad\blacksquare$}}

\newcommand{\Op}{{\mathcal O}}
\newcommand{\Br}{{\mathcal B}_r}
\newcommand{\Bc}{{\mathcal B}_c}
\newcommand{\acn}{{\mathcal A}^n_c}
\newcommand{\arn}{{\mathcal A}^n_r}
\newcommand{\ac}{{\mathcal A}_c}
\newcommand{\ar}{{\mathcal A}_r}

\newcommand{\ftilde}{{\tilde f}}
\newcommand{\Ftilde}{{\tilde F}}

\newcommand{\supp}{{\rm supp}}

\newcommand{\Rbar}{\overline{\R}}

\newcommand{\D}{{\mathcal D}(\R)}
\newcommand{\Dp}{{\mathcal D}'(\R)}

\newcommand{\nbvl}{{\mathcal NBV}_{\!\lambda}}

\newcommand{\ebv}{{\mathcal EBV}}
\newcommand{\bv}{{\mathcal BV}}
\newcommand{\intinf}{\int^\infty_{-\infty}}

\newcommand{\N}{{\mathbb N}}
\newcommand{\R}{{\mathbb R}}
\newcommand{\C}{{\mathbb C}}

\newcommand{\Sc}{{\mathcal S}}
\newcommand{\fn}{\!:\!}
\newcommand{\lsum}{\sum\limits}

\newcommand{\lint}{\int\limits}

\DeclareMathOperator{\nbvone}{{\mathcal NBV}_{\!1}}
\DeclareMathOperator{\nbvzero}{{\mathcal NBV}_{\!0}}
\DeclareMathOperator{\esssup}{ess\,sup}

\providecommand{\norm}[1]{\lVert#1\rVert}
\newtheorem{theorem}{Theorem}
\newtheorem{lemma}[theorem]{Lemma}
\newtheorem{prop}[theorem]{Proposition}

\newtheorem{defn}[theorem]{Definition}

\begin{document}
\hspace{-2cm}
\raisebox{12ex}[1ex]{\fbox{{\footnotesize
Preprint
October 12, 2011.\quad
To appear in {\it Czechoslovak Mathematical Journal}
}}}

\begin{center}
{\large\bf Integrals and Banach spaces for finite order distributions}
\vskip.25in
Erik Talvila\footnote{Supported by the
Natural Sciences and Engineering Research Council of Canada.
}\\ [2mm]
{\footnotesize
Department of Mathematics and Statistics \\
University of the Fraser Valley\\
Abbotsford, BC Canada V2S 7M8\\
Erik.Talvila@ufv.ca}
\end{center}

{\footnotesize
\noindent
{\bf Abstract.} 
Let $\Bc$ denote the real-valued functions continuous on the
extended real line and vanishing at $-\infty$.  Let $\Br$ denote
the functions that are left continuous, have a right limit
at each point and vanish at $-\infty$.  Define $\acn$ to be the
space of tempered distributions that are the $n$th distributional
derivative of a unique function in $\Bc$.  Similarly with
$\arn$ from $\Br$.  A type of integral is defined on distributions
in $\acn$ and $\arn$.  The multipliers are iterated integrals of
functions of bounded variation.  For each $n\in\N$, the spaces
$\acn$ and $\arn$ are Banach spaces, Banach lattices and Banach
algebras isometrically isomorphic to $\Bc$ and $\Br$, respectively.
Under the ordering in this lattice, if a distribution is integrable
then its absolute value is integrable.
The dual space is isometrically isomorphic to the functions of 
bounded variation.  The space $\ac^1$ is
the completion of the $L^1$ functions in the Alexiewicz norm.  The
space $\ar^1$ contains all finite signed Borel measures.  Many of
the usual properties of integrals hold: H\"older inequality,
second mean value theorem, continuity in norm, linear change of variables,
a convergence theorem.
\\
{\bf 2010 subject classification:} 26A39, 46B42, 46E15, 46F10, 46G12,
46J10
}\\

\section{Introduction}\label{sectionintroduction}
An integral means different things to different people: a Riemann
sum, approximation
by simple functions, inversion of the derivative of an absolutely
continuous function,
a type of linear functional.  In this paper
we will define an integration process on Schwartz distributions
of finite order
by inverting the $n$th order distributional derivative of regulated or
continuous
functions.  An
important part of the definition is that we will obtain a linear
functional that acts on iterated integrals of functions of bounded
variation.

We denote by $\Bc$ the continuous functions on the extended real line that
vanish at $-\infty$.
Then $\acn$ is the set of distributions that are the $n$th
distributional derivative of a unique function in $\Bc$.  The case $n=1$
was studied in \cite{talviladenjoy}.  The space $\ac^1$ is the
completion of the space of $L^1$ functions and the completion of the
space of Henstock--Kurzweil integrable functions in the Alexiewicz norm.
If $f\fn\R\to\R$ is a Henstock--Kurzweil integrable function then
the Alexiewicz norm is $\norm{f}=\sup_{x\in\R}|\int_{-\infty}^xf(t)\,dt|$.
If $F\in\Bc$ and $f=F'$ is
its distributional derivative then the {\it continuous primitive
integral} of $f$ is $\int_a^bf=F(b)-F(a)$.  In this paper we define
$\intinf fh$ for $f\in\acn$ and $h$ an $n$-fold iterated integral of
a function of bounded variation (Definition~\ref{defnint}).  
A function is called {\it regulated} if it has
a left limit and a right limit at each point.  We denote by $\Br$ the 
left continuous regulated functions that
vanish at $-\infty$.  Then $\arn$ is the set of distributions that are the $n$th
distributional derivative of a function in $\Br$.  The case $n=1$
was studied in \cite{talvilaregulated}.  The space $\ar^1$ contains
$\ac^1$ as well as all finite signed Borel measures. An integral in $\arn$ is
defined as for $\acn$.  Under the uniform norm, $\Bc$ and $\Br$ are 
Banach spaces.  The $n$th order distributional derivative provides a linear
isometric
bijection between $\Bc$ and $\acn$ and between $\Br$ and $\arn$.  For each
$n\in\N$,
the spaces of distributions $\acn$ and $\arn$ are then Banach spaces
that are isometrically isomorphic to $\Bc$ and $\Br$, respectively.
If $f\in\arn$ with primitive $F\in\Br$ then its norm is $\norm{f}_{a,n}=
\norm{F}_\infty$.  The spaces $\Bc$ and $\acn$ are separable while the
spaces $\Br$ and $\arn$ are not separable.  

Below we define our notation for distributions.  

The main definitions are given in Section~\ref{sectiondefn}.  For $f\in\arn$,
the
integral $\intinf fh$ is defined by reducing to $\intinf F'h^{(n-1)}$.
With $h^{(n-1)}$ a function of bounded variation, this reduces to an integral
in $\ar^1$.  This is evaluated using Henstock--Stieltjes integrals.
See Definition~\ref{defnint}.
  It is
shown that distributions in $\arn$ are tempered and of order at most
$n$. There is translation invariance and for $\acn$ there is continuity in norm.
The
multipliers for these integrals are iterated integrals of functions
of bounded variation.  Various properties of these functions are
proved here.

Examples and further properties of the integral are given in
Section~\ref{sectionexamples}.  If $\delta$ is the Dirac distribution
then $\delta^{(m)}\in\ar^{m+1}$ for each $m\geq 0$.  It is shown that
each distributional derivative of a finite signed Borel measure is
in some $\arn$ space.  A version of the second mean value
theorem is established and a linear change of variables theorem is proved.  
Relationships amongst $\acn$ and $\ar^m$ are investigated.

In Section~\ref{sectionholder} a type of H\"older inequality is
established.  For $f\in\arn$ with primitive $F\in\Br$ we have
$|\intinf fh|\leq \norm{F}_\infty \norm{g}_\bv$, where $g$ is
of bounded variation and $g=h^{(n-1)}$. This leads to a convergence
theorem when a sequence $(f_n)\subset\arn$ converges in norm.
Using the isometries $\acn\leftrightarrow\Bc$ and $\arn\leftrightarrow\Br$
it is shown that the dual space of $\acn$ is the functions of normalised
bounded variation and the dual space of $\arn$ is the functions of bounded
variation.

Under pointwise operations $\Br$ is a Banach
lattice.  If $F_1,F_2\in\Br$ then $F_1\leq F_2$ means $F_1(x)\leq F_2(x)$ 
for all $x\in\R$.
This Banach lattice structure is inherited by $\acn$ and $\arn$.  If
$f_1=F_1^{(n)}, f_2=F_2^{(n)}\in\arn$ then $f_1\preceq f_2$ if and only
if $F_1\leq F_2$ in $\Br$.  Elementary lattice properties are proved in
Section~\ref{sectionlattice},
including the fact that $\acn$ and $\arn$ are abstract $M$ spaces.
Under this lattice ordering, the integrals introduced are absolute
in the sense that if $f$ is integrable then the absolute value
of $f$ in this ordering is integrable.
This is the case even though functions in $\arn$ may have conditionally
convergent Henstock--Kurzweil integrals.

Banach algebras are considered in Section~\ref{sectionalgebra}.
Under pointwise operations $\Br$ is a Banach
algebra.  If $F_1,F_2\in\Br$ then $(F_1F_2)(x)=F_1(x)F_2(x)$.
By the isomorphism, $\acn$ and $\arn$ are also Banach algebras.
If
$f_1=F_1^{(n)}, f_2=F_2^{(n)}\in\arn$ then $f_1f_2=D^n(F_1F_2)$.
For
complex-valued distributions they are  $C^\ast$-algebras. 
Under this multiplication $\delta^{(n)}\delta^{(n)}=\delta^{(n)}$ for
each $n\geq 0$.
J.F.~Colombeau, E.E.~Rosinger and others have embedded spaces of
distributions in various algebras.  For example, see \cite{oberguggenberger}.
The Banach algebra we construct here seems to be unrelated.

The starting point in this paper are the spaces $\Bc$ and $\Br$.
Basic properties are established in Section~\ref{sectiondefn}.  The
following lemma is used repeatedly to carry over Banach space,
Banach lattice and Banach algebra properties to $\acn$ and $\arn$.
\begin{lemma}\label{lemma}
Let $A$ be a set.  Let $B$ be a vector space over field $\R$.  
Let $x,y\in A$; $a\in \R$.\\
(a) Suppose there is a bijection $\Phi\fn B\to A$.
Define $x+y=\Phi(\Phi^{-1}(x)+\Phi^{-1}(y))$ and $ax=\Phi(a\Phi^{-1}(x))$.
Then $A$ is vector space isomorphic to $B$ and $\Phi$ is linear such that
$\Phi\circ\Phi^{-1}=i_A$ and $\Phi^{-1}\circ\Phi=i_B$.\\
(b) Suppose $B$ is a Banach space.  Define
$\norm{x}_A=\norm{\Phi^{-1}(x)}_B$.  Then $A$ is a Banach space isometrically
isomorphic to $B$.  If $B$ is separable so is $A$.\\
(c)  Suppose $B$ is a Banach lattice.  Define
$x\preceq y$ in $A$ if and only if $\Phi^{-1}(x)\preceq\Phi^{-1}(y)$ in $B$.  Then
$A$ and $B$ are isometrically isomorphic Banach lattices.\\
(d)   Suppose $B$ is a Banach algebra.  Define
$xy=\Phi(\Phi^{-1}(x)\Phi^{-1}(y))$.  Then 
$A$ and $B$ are isomorphic Banach algebras.
\end{lemma}
The proof is elementary.   A related result is that if $\Phi$ is a
surjective isometry between two normed linear spaces then $\Phi$ must
be linear.  This is the
Mazur--Ulam theorem.  For example, see \cite{fleming}.

The {\it test functions} are $\D=C^\infty_c(\R)$,  i.e., the smooth functions
with compact support.  The {\it support} of  a  function $\phi$ is the closure
of  the  set on  which $\phi$ does not  vanish.  Denote this as $\supp(\phi)$.
There is a notion of  continuity in $\D$.  If
$\{\phi_n\}\subset \D$ then $\phi_n\to\phi\in\D$ if there is a compact
set $K\subset\R$  such  that for  all  $n\in\N$, $\supp(\phi_n)\subset K$,
and for each integer $m\geq 0$, $\phi^{(m)}_n\to \phi^{(m)}$ uniformly
on $K$ as $n\to\infty$.  The {\it distributions} are the  continuous linear
functionals on $\D$,  denoted $\Dp$.  If $T\in\Dp$ then $T\fn\D\to\R$
and we write $\langle T,\phi\rangle\in\R$ for $\phi\in\D$.  If
$\phi_n\to\phi$ in $\D$ then $\langle T,\phi_n\rangle\to\langle T,\phi\rangle$
in $\R$.  And, for all $a_1, a_2\in\R$ and all $\phi,\psi\in\D$, $\langle T,
a_1\phi+a_2\psi\rangle =a_1\langle T,\phi\rangle+a_2\langle T,\psi\rangle$.  If
$f\in L^p_{{\rm loc}}$  for some $1\leq p\leq\infty$ then $\langle T_f,\phi\rangle=\int_{-\infty}^\infty
f(x)\phi(x)\,dx$ defines a distribution.  For a locally integrable function we
will often drop the distinction between $f$ and $T_f$. The  
differentiation  formula $\langle D^nT,\phi\rangle=\langle T^{(n)},
\phi\rangle=(-1)^n\langle T,\phi^{(n)}\rangle$ ensures that all 
distributions have
derivatives of all  orders which are themselves  distributions.  This is
known as the {\it distributional derivative} or {\it weak derivative}.
We  will  usually denote distributional
derivatives by $D^nF$, $F^{(n)}$ or $F'$ and  pointwise derivatives by $F^{(n)}(t)$
or $F'(t)$.
For $T\in\Dp$ and $t\in\R$  the {\it translation} $\tau_t$ is defined
by $\langle\tau_tT, \phi\rangle=\langle T, \tau_{-t}\phi\rangle$
where $\tau_t\phi(x)=\phi(x-t)$ for
$\phi\in\D$.  If there is an integer $N\geq 0$ such that for each compact
set $K\subset\R$ there is a real number $C\geq 0$ so that 
$|\langle T,\phi\rangle|\leq C\sum_{0}^N\norm{\phi^{(n)}}_\infty$ for
all $\phi\in\D$ with support in $K$, then distribution $T$ is said to 
be of {\it finite order}.
The least such $N$ is the {\it order} of $T$.
If $\mu$ is a finite signed Borel measure
then $\langle T_\mu,\phi\rangle=\intinf\phi(x)\,d\mu(x)$ defines
$T_\mu\in\Dp$ as a distribution of order $0$.
Most of the results on distributions we
use
can  be found in \cite{folland}, \cite{friedlander} or \cite{zemanian}.

Several authors have proposed various schemes for integrating
distributions.  L.~Schwartz \cite{schwartz} considered the integral of $T\in\Dp$
as the linear functional $\langle T,1\rangle$, whenever this exists.
As will be seen in the next section, we generalise 
Schwartz's definition so that the integrable distributions are
continuous
linear functionals on iterated integrals of functions of bounded
variation.
J.C.~Burkill \cite{burkill} has sketched out a method of integrating
distributions using higher order Stieltjes integrals.  
A.M.~Russell \cite{russell}
and A.G.~Das \cite{das} with coauthors have also used higher Stieltjes 
integrals.  
J.~Mikusi\'nski, J.A.~Musielak and R.~Sikorski
have used convolutions to define a type of integral for distributions.
See
\cite{jmikusinski}, \cite{musielak}, \cite{sikorski}.

The extended real line is denoted $\Rbar=[-\infty,\infty]$.
The space
$C(\Rbar)$ consists of the continuous functions $F\fn\Rbar\to\R$.
A function is in $C(\Rbar)$ if it is continuous at each point in 
$\R$ and if
$F(\infty)=\lim_{x\to\infty}F(x)\in\R$ and $F(-\infty)=\lim_{x\to-\infty}F(x)
\in\R$.  This two-point compactification makes $\Rbar$ into a compact
Hausdorff space.  A topological
base for $\Rbar$ consists of the usual open intervals $(a,b)$ with
$-\infty\leq
a<b\leq\infty$, as well as $[-\infty,a)$ with $-\infty<a\leq\infty$,
and $(a,\infty]$ with $-\infty\leq a<\infty$.   In this paper the
word {\it compact} will always refer to the usual topology on $\R$.

A function
$F\fn\Rbar\to\R$ is {\it regulated} on $\Rbar$ if it has a left and right limit
at each point of $\R$ and real limits at $\pm\infty$, i.e., for
each $x\in\R$ the limits $F(x-)=\lim_{y\to x^-}F(y)$ and 
$F(x+)=\lim_{y\to x^+}F(y)$ exist as real numbers and 
$\lim_{y\to-\infty}F(y)$ and $\lim_{y\to\infty}F(y)$ exist as
real numbers.  We will use the following normalisations for regulated
functions.  If $F$ is regulated and $0\leq\lambda\leq 1$ then
$F_\lambda(x)=(1-\lambda) F(x-)+\lambda F(x+)$ 
for all $x\in\R$.  The functions $F$ and $F_\lambda$ will then
differ on a countable set.  Note that $F_0$
is left continuous and 
$F_1$ is right continuous.  If $F$ is continuous then
all normalisations are equal to $F$.  The {\it Heaviside step function} 
has the left and right continuous normalisations
$H_0=\chi_{(0,\infty]}$ and $H_1=\chi_{[0,\infty]}$.  Unless otherwise
stated, all regulated functions will satisfy $F(-\infty)=
\lim_{y\to-\infty}F(y)$ and $F(\infty)=\lim_{y\to\infty}F(y)$.
For more on regulated functions, see \cite{frankova}.

\section{Banach spaces and integrals}\label{sectiondefn}
A space of primitives is the
regulated functions that vanish at $-\infty$.
Each such regulated function is differentiated
$n$ times with the distributional derivative to yield a sequence of 
Banach spaces of integrable distributions, each being isometrically
isomorphic to the space of primitives.
A second space of primitives is the set of functions in $C(\Rbar)$ that
vanish at $-\infty$.  These are also differentiated
$n$ times to give a sequence of Banach spaces.
By integrating functions of bounded
variation $n$ times we find the corresponding set of multipliers.
Distinction is made between functions of bounded variation, normalised
bounded variation and essential bounded variation.

We take as our set of primitives
$\Br$.  This consists of the functions $F\fn\Rbar\to\R$ that are regulated
on $\Rbar$
such that $F(-\infty)=\lim_{x\to-\infty}F(x)=0$, $F(x)=F(x-)$ for all $x\in\R$ and
$F(\infty)=\lim_{y\to\infty}F(y)$.  
Hence, they
are left continuous on $(-\infty,\infty]$, vanish at $-\infty$ and 
equal their limits at infinity.
Under pointwise operations and
the uniform norm, $\norm{F}_\infty=\sup_{x\in\R}|F(x)|$, 
$\Br$ is a Banach space.  It is not separable.  See
\cite[Theorem~2]{talvilaregulated} where various properties of $\Br$ are
proved.
A second set of primitives is
the subspace $\Bc=\{F\in C(\Rbar)\mid F(-\infty)=0\}$.  This is then a
separable Banach
space with norm $\norm{F}_\infty=\sup_{x\in\R}|F(x)|=\max_{x\in\Rbar}|F(x)|$.
The separability of $\Bc$ follows from the compactness of $\Rbar$.
See \cite[Exercise V.7.12]{dunfordschwartz}.

For each $n\in\N$ define 
$\arn=\{f\in\Dp\mid f=F^{(n)} \text{ for some } F\in \Br\}$,
i.e., $\langle f,\phi\rangle=\langle F^{(n)},\phi\rangle=(-1)^n\langle F,\phi^{(n)}\rangle
=(-1)^n\intinf F(x)\phi^{(n)}(x)\,dx$
for each $\phi\in\D$.  This last integral is a Riemann integral with a
compactly supported integrand.  From 
this definition we see that elements of $\arn$ are distributions.
And, define $\acn=\{f\in\Dp\mid f=F^{(n)} \text{ for some } F\in \Bc\}$.

\begin{theorem}[Uniqueness]\label{thmunique}
For each $f\in\arn$ there is a unique function $F\in\Br$ such that $F^{(n)}=f$.
For each $f\in\acn$ there is a unique function $F\in\Bc$ such that $F^{(n)}=f$.
\end{theorem}

\bigskip
\noindent
{\bf Proof:}
Suppose $F^{(n)}=G^{(n)}$ for some $F,G\in\Br$ then let $P=F-G$.
Thus, $P\in\Br$, $P^{(n)}=0$ and $P$ is a polynomial of degree at most $n-1$.
This follows from the fact that the polynomials are Fourier transforms
of linear combinations of the Dirac distribution and its derivatives
\cite[Exercise~9.25]{folland}.  The only polynomial in $\Br$ is $0$.
Similarly when $F,G\in\Bc$.
\qed

If $f\in\arn$
we can then speak of the unique element $F\in\Br$ such that $F^{(n)}=f$ as 
the {\it primitive} of
$f$.  Here it is essential that $F$ be left continuous rather than just
regulated.  The mapping $\Phi\fn\Br\to\arn$ given by $\Phi(F)=F^{(n)}$ is a linear bijection.
It is surjective by the definition of $\arn$.  It is injective by Theorem~\ref{thmunique}.
It follows from Lemma~\ref{lemma} that $\arn$ is a linear space.  The
norm inherited from $\Br$ makes $\arn$ into a Banach space.  Similarly
for $\acn$.
We call this the Alexiewicz norm and denote it $\norm{\cdot}_{a,n}$.
(See \cite{alexiewicz}.)
The Alexiewicz norm is translation invariant and we have continuity
in norm in $\acn$ but not in $\arn$.
And, $C^\infty(\R)$ is dense in $\acn$ but not in $\arn$.
\begin{theorem}
Let $n\in\N$ and
let $f, f_1, f_2\in\arn$ with respective primitives $F, F_1, F_2\in\Br$.
Let $a_1,a_2\in\R$.  Let $\phi\in\D$.
(a) With
operations 
given
by 
\begin{eqnarray*}
\langle a_1f_1+a_2f_2,\phi\rangle & = & a_1\langle f_1,\phi\rangle
+a_2\langle f_2,\phi\rangle\\
 & = &  (-1)^n(a_1\langle F_1,\phi^{(n)}\rangle+
a_2\langle F_2,
\phi^{(n)}\rangle)\\
 & = & (-1)^n(a_1\intinf F_1(x)\phi^{(n)}(x)\,dx
+a_2\intinf F_2(x)\phi^{(n)}(x)\,dx)
\end{eqnarray*}
$\arn$ is a vector space.  And, $\acn$ is a subspace of $\arn$.
(b) A norm on $\arn$ is defined by $\norm{f}_{a,n}=\norm{F}_\infty$.
This makes $\arn$ into a Banach space that is not separable.
Each of the spaces $\arn$ is isometrically isomorphic to $\Br$.
Each of the spaces $\acn$ is a separable Banach space isometrically
isomorphic to $\Bc$.
(c) Each distribution in $\arn$ is tempered, of order at most $n$.
(d) Let $\Op\fn\Dp\to\Dp$ be an operator that commutes
with the derivative, $(\Op T)'=OT'$ for all $T\in\Dp$.  Then
$\langle \Op f,\phi\rangle=(-1)^n\langle \Op F, \phi\rangle$.
(e) Let $t\in\R$ and let  $T\in\Dp$. Then $T\in\arn$ if and only if 
$\tau_tT\in\arn$.  Similarly for $\acn$.  
If $f\in\arn$ then $\norm{\tau_tf}_{a,n}=\norm{f}_{a,n}$.
(f) For each $f\in\acn$ it follows that $\lim_{t\to 0}\norm{f-\tau_tf}_{a,n}=0$.
(g) $C^\infty(\R)$ is dense in $\acn$ but not in $\arn$.
\end{theorem}

\bigskip
\noindent
{\bf Proof:}
(a) and (b) These follow from Lemma~\ref{lemma}. 

(c) Let $K\subset\R$ be compact.  Suppose $\phi\in\D$
with $\supp(\phi)\subset K$.  Denote the Lebesgue measure of $K$
by $|K|$.  Let $f\in\arn$ with primitive $F\in\Br$.
Then
$$
|\langle f,\phi\rangle|=\left|\intinf F(x)\phi^{(n)}(x)\,dx\right|\leq
\norm{F}_\infty|K|\norm{\phi^{(n)}}_\infty.$$
This shows $f$ is tempered and of order at most $n$.  See \cite{friedlander}
for the definition of tempered.

(d) It follows from associativity that $\Op T^{(n)}=(\Op T)^{(n)}$.  Then
$\langle \Op f,\phi\rangle=\langle\Op F^{(n)},\phi\rangle=\langle
(\Op F)^{(n)},\phi\rangle=(-1)^{n}\langle \Op F,\phi^{(n)}\rangle$.

(e)  Let $T\in\Dp$.  Then $\langle\tau_tT',\phi\rangle=\langle T',\tau_{-t}\phi\rangle
=-\langle T, (\tau_{-t}\phi)'\rangle$.
And,
$\langle (\tau_tT)',\phi\rangle=-\langle\tau_tT,\phi'\rangle=-\langle T,
\tau_{-t}\phi'\rangle$.  For $x\in\R$, $\tau_{-t}\phi'(x)=(\tau_{-t}\phi)'(x)
=\partial\phi(x+t)/\partial x$.  It follows that
$\tau_t$ commutes with derivatives.  From the proof of (d), if $f\in\arn$ then
$\langle\tau_tf,\phi\rangle   =  \langle\tau_t F^{(n)},\phi\rangle=
   \langle(\tau_t F)^{(n)}, \phi\rangle$.
Therefore, $\tau_tf$ is the $n$th derivative of $(\tau_tF)\in\Br$,
so $\tau_t f\in\arn$.  If $T\in\Dp$ such that $\tau_tT\in\arn$ then
write $T=\tau_{-t}(\tau_tT)$ to show $T\in\arn$.
Similarly for $\acn$.
For $f\in\arn$, $\norm{\tau_tf}_{a,n}=\norm{\tau_tF}_\infty=\norm{F}_\infty=
\norm{f}_{a,n}$.

(f) We have $\norm{f-\tau_tf}_{a,n}=\norm{
F^{(n)}-(\tau_tF)^{(n)}}_{a,n}=\norm{F-\tau_tF}_\infty\to 0$ as $t\to 0$.

(g) Let $\Phi_y(x)=(y/\pi)(x^2+y^2)^{-1}$ be the half plane Poisson
kernel.  Define the convolution $G_y=F\ast\Phi_y$.  
Since $F$ is continuous on $\Rbar$, it is known that $\norm{G_y
-F}_\infty\to0$ as $y\to 0^+$.  See, for example, \cite{axler}.
Note that $G_y\in C^\infty(\R)$.  By
dominated convergence and the fact that $\intinf \Phi_y(x)\,dx=1$ we
see that $\lim_{x\to\infty}G_y(x)=F(\infty)$ and
$\lim_{x\to-\infty}G_y(x)=F(-\infty)=0$.  Hence, $G_y\in\Bc$ for 
each $y>0$.  The density of $C^\infty(\R)$ in $\acn$ now follows.
\qed

Besides translations, other
examples of operators commuting with the derivative are
linear
combinations of differential operators with coefficients independent
of the differentiation variable.

It was shown in \cite[Proposition~3.3]{talvilaconvolution} that $L^1$, 
and hence the space
of Henstock--Kurzweil integrable functions, is dense in $\ac^1$.
Note that $C^\infty(\R)$ is not dense in $\arn$.
The Heaviside step function
$H_0$ is in $\Br$.  For each $\psi\in C(\R)$ we have $\norm{H_0-\psi}_\infty
\geq 1/2$.  Therefore, $C^\infty(\R)$ is not dense in $\arn$.
And, $\D=C^\infty_c(\R)$ is not dense in $\acn$.  Define
$F\in\Bc$ by
$F(x)=\pi/2+\arctan(x)$.  Then for each $\phi\in\D$ we have
$\norm{F-\phi}_\infty\geq \pi$.

We do not have continuity in norm in $\arn$.  For example, consider
$\norm{H_0^{(n)}-\tau_tH_0^{(n)}}_{a,n}=\norm{H_0-\tau_tH_0}_\infty=1$
for all $t\not=0$.  

The {\it variation} of a function
$g\fn\Rbar\to\R$ is the supremum of $\sum|g(x_i)-g(y_i)|$, 
taken over all disjoint
intervals $(x_i,y_i)\subset\R$.  This is denoted $V g$.  The 
functions of {\it bounded variation} are
$\bv=\{g\fn\Rbar\to\R\mid Vg<\infty\}$.
Functions of bounded variation are the difference of two increasing
functions and thus are regulated on $\Rbar$.
Under usual pointwise operations $\bv$ is a Banach space under the norm
$\norm{g}_{\bv}=\norm{g}_\infty+Vg$.  For each $-\infty\leq a\leq \infty$,
an equivalent norm is $|g(a)|+Vg$.  See, for
example, \cite{dunfordschwartz}, \cite{folland} and \cite{kannan} 
for properties of $\bv$ functions.

The following spaces will serve as multipliers for $\arn$.  
Each space $\Ibvn$ is defined
inductively. 
\begin{defn}
Define $\Ibv^0=\bv$.  Suppose $\Ibv^{n-1}$ is known for $n\in\N$.
Define $\Ibvn=\{h\fn\R\to\R\mid
h(x)=\int_0^xq(t)\,dt \text{ for some } q\in\Ibv^{n-1}\}$.
\end{defn}
Hence, a function $h\in\Ibvn$ is an $n$-fold iterated integral.
\begin{defn}\label{defnI^n}
Let $q\in L^1_{loc}$.  Define
$I^0[q](x) =q(x)$. For $n\in\N$  define
$$
I^n[q](x)=\lint_{x_n=0}^x\cdots\lint_{x_i=0}^{x_{i+1}}\cdots\lint
_{x_1=0}^{x_2}q(x_1)\,dx_1\cdots dx_{i}\cdots dx_n.
$$
\end{defn}

\begin{prop}\label{propIbvnint}
If $h\in\Ibvn$ for $n\in\N$ then there is a function $g\in\bv$ such that for 
all $x\in\R$, 
$$
h(x)=I^n[g](x)=
\frac{1}{(n-1)!}\int_0^x(x-s)^{n-1}g(s)\,ds.
$$
\end{prop}
The proof follows from induction and the Fubini--Tonelli theorem.

The function $g$ is not unique since there are functions of bounded
variation that differ only on a countable set.  Imposing a normalisation
on $\bv$ makes $g$ unique.  Fix $0\leq \lambda\leq 1$.  
Functions $g$ and $g_\lambda$ differ on a set
that is countable and $Vg_\lambda=\inf Vh$ where the infimum is
taken over all $h\in\bv$ such that $g_\lambda=h$ almost everywhere.
The value of $0\leq \lambda\leq 1$ does not affect the value of
$Vg_\lambda$.
The functions of {\it normalised bounded variation}
are then 
$\nbvl=\{g_\lambda\mid g\in\bv\}$.  It is easy to see that for
each $0\leq \lambda\leq 1$ there is a
unique $g\in\nbvl$ such that if $h\in\Ibvn$ then $h=I^n[g]$.
Clearly, $\nbvl\subsetneq\bv$.

Note that if $g\in\bv$ then it is bounded so the function
$x\mapsto \int_0^xg(t)\,dt$ is Lipschitz continuous and $I^n[g]\in C^{n-1}(\R)$.
The same applies if $g\in\ebv$ (see below).
If $h\in\Ibvn$ then $h^{(m)}(0)=0$ for all $0\leq m\leq n-1$ 
and $h(x)=O(x^n)$ as $|x|\to\infty$.
 
Now we can define integrals on $\arn$.
A distribution $f\in\ar^1$ is the
distributional derivative of a unique function $F\in\Br$.  Its
{\it regulated primitive integrals} are
\begin{eqnarray}
\int_{(a,b)}f & = & \int_{a+}^{b-}f=F(b-)-F(a+) =F(b)-F(a+)\label{regulatedintdefn1}\\
\int_{(a,b]}f & = & \int_{a+}^{b+}f=F(b+)-F(a+)\label{regulatedintdefn2}\\
\int_{[a,b)}f & = & \int_{a-}^{b-}f=F(b-)-F(a-)=F(b)-F(a)\label{regulatedintdefn3}\\
\int_{[a,b]}f & = & \int_{a-}^{b+}f=F(b+)-F(a-)
=F(b+)-F(a)\label{regulatedintdefn4}
\end{eqnarray}
for all $-\infty< a<b< \infty$.  We also have $\intinf f=F(\infty)$ with
similar definitions for semi-infinite intervals.  And,
$\int_{\{a\}}f=F(a+)-F(a-)$.  This integral was described in
detail in \cite{talvilaregulated}.
The multipliers are the functions of bounded variation. If $g\in\bv$
then 
\begin{eqnarray}
\intinf fg & = & \intinf g(x)\,dF(x)=F(\infty)g(\infty)-\intinf F(x)\,dg(x)\notag\\
 &  & \quad  -\sum_{n\in\N}\left[
F(c_n)-F(c_n+)\right]\left[g(c_n)-g(c_n+)\right].\label{ar1defn}
\end{eqnarray}
The sum is over all $c_n\in\R$ at which $F$ and $g$ are not both
 right continuous.
The integrals $\intinf g(x)\,dF(x)$ and $\intinf F(x)\,dg(x)$ are
Henstock--Stieltjes integrals.  They are known to exist when one of
$F$ and $g$ is regulated and one is of bounded variation.  They are defined
by using tagged partitions of $\Rbar$.
If $F$ is regulated but not required to be left continuous then an
additional term containing coincident jump discontinuities of $F$ and
$g$ from the left must be added.  See \cite[p.~199]{mcleod}.

If $f\in\ac^1$ then the four integrals 
\eqref{regulatedintdefn1}-\eqref{regulatedintdefn4}
all give $\int_a^bf=F(b)-F(a)$ and the sum in \eqref{ar1defn} vanishes.
In this case, $\intinf F(x)\,dg(x)$ is a Riemann--Stieltjes integral over an
unbounded domain.
It can also be defined by taking limits of 
Riemann--Stieltjes integrals
over finite subintervals.
See \cite[p.~187]{mcleod} and \cite{talvilaraelimits} for details.  

\begin{defn}\label{defnint}
Let $n\in\N$ and let $0\leq\lambda\leq 1$.  For $f\in\arn$ let $F$ be its primitive in $\Br$.
For $h\in\Ibv^{n-1}$ such that $h=I^{n-1}[g]$ for $g\in\nbvl$, define
the regulated primitive integral of $f$ with respect to $h$ as
\begin{eqnarray}
\intinf fh & = & \intinf F^{(n)}h=(-1)^{n-1}\intinf F'h^{(n-1)}
\label{intdefn1}\\
 & = & (-1)^{n-1}F(\infty)g(\infty)-(-1)^{n-1}\intinf F(x)\,dg(x)
\label{intdefn2}\\
 & & \quad
-(-1)^{n-1}\sum_{n\in\N}\left[
F(c_n)-F(c_n+)\right]\left[g(c_n)-g(c_n+)\right]
\label{intdefn3}\\
 & = & (-1)^{n-1}\intinf h^{(n-1)}(x)\,dF(x).\label{intdefn4}
\end{eqnarray}
\end{defn}

If $f\in\acn$, or if $f\in\arn$ and $g\in\nbvone$ (i.e. right
continuous), then the sum in 
\eqref{intdefn3} vanishes.

To distinguish them from the regulated primitive integral, we will
always explicitly show the integration variable and differential
in Lebesgue, Henstock--Stieltjes
and Riemann integrals.
It is shown in \cite{talviladenjoy} (following Definition~6) 
that if $g_1, g_2\in\bv$ differ
on a countable set and $F\in C(\Rbar)$ then $\intinf F(x)\,dg_1(x)=
\intinf F(x)\,dg_2(x)$.  Hence, if $f\in\acn$ it makes no difference in 
Definition~\ref{defnint} if we use $g\in\nbvl$ or any function of
bounded variation that differs from $g$ on a countable set.
\begin{prop}
Let $n\in\N$.  For $f\in\acn$ let $F$ be its primitive in $\Bc$.
Let $h\in\Ibv^{n-1}$ and let $g\in\bv$ such that $h=I^{n-1}[g]$.  Then
$\intinf fh= (-1)^{n-1}F(\infty)g(\infty)-(-1)^{n-1}\intinf F(x)\,dg(x)$.
\end{prop}

The integral in \eqref{intdefn2} takes different forms when
$g$ is the primitive for different types of integrals.
The set of primitives for $L^1$ functions is $AC(\Rbar):=AC(\R)\cap\bv$,
where $AC(\R)$ are the functions that are absolutely continuous on
each compact set in $\R$.  If $f$ is a measurable function on the real
line then $f\in L^1$ if and only if there is a function $F\in AC(\Rbar)$
such that $f(x)=F'(x)$ for almost all $x\in\R$.
\begin{prop}\label{propLebesgueRiemann}
Let $F\in \Br$.
(a) If $g\in AC(\Rbar)$ then
$\intinf F(x)\,dg(x)=
\intinf F(x)g'(x)\,dx$.  This last is a Lebesgue integral.
(b)  If $g\in C^1(\R)\cap C(\Rbar)$ then
$\intinf F(x)\,dg(x)=
\intinf F(x)g'(x)\,dx$.  This last is an improper Riemann integral.
\end{prop}
This follows from the form the fundamental theorem of calculus
takes for each integral.  For example, see \cite[p.~74]{leevyborny}.
The primitives for the Riemann integral are the functions
of bounded slope variation.  See \cite{thomson}.  This set of primitives
properly contains the Lipschitz functions and is a proper subset of  
$C^1(\R)$.

Let $g\fn\R\to\R$.  Then $g$ is of {\it essential bounded
variation} if its distributional derivative is a signed Radon
measure, i.e., there is a unique signed Radon measure $\mu_g$ such
that $\langle g',\phi\rangle = -\langle g,\phi'\rangle
=-\intinf g(x)\phi'(x)\,dx=\intinf \phi(x)\,d\mu_g(x)$ for all
$\phi\in\D$.  Denote the functions of essential bounded
variation by $\ebv$.  If $g\in\ebv$ then its {\it essential
variation} is $EVg=\sup_\phi\intinf g(x)\phi'(x)\,dx$, the
supremum being taken over all functions $\phi\in\D$ such
that $\norm{\phi}_\infty\leq 1$.  Changing a function on a set of
measure zero does not affect its essential variation. 
The elements of $\ebv$ are equivalence classes of functions 
that are equal almost everywhere.  For each $0\leq\lambda\leq 1$
there is exactly one function from $\nbvl$ in each equivalence class.
If $g\in\ebv$ then
there is
exactly one function
$h\in\nbvl$ such that $EVg=Vh$.  
Hence, $\ebv$ and $\nbvl$ are isometrically isomorphic. 
The space
$\ebv$ is a Banach space under norm $\norm{g}_{{\mathcal EBV}}= \esssup|g|
+EVg=\norm{h}_{\bv}$.  Also, $EVg=|\mu_g|(\R)=\inf_h Vh$
where the infimum is taken over all $h\in\bv$ such that
$g=h$ almost everywhere.
It is shown in \cite[Corollary~15]{talvilaconvolution} that
$\intinf F(x)\,dg(x)=\intinf F(x)\,d\mu_g(x)$ for all $F\in C(\Rbar)$.
A limiting process is used in \cite[Theorem~8]{talviladenjoy} to
define $\intinf F'g$ for $F\in\Bc$ and $g\in\ebv$.  Hence,
in Definition~\ref{defnint} we can use $g\in\ebv$ when $f\in\acn$.
All of the results in this paper can be rewritten under this
assumption.  If $h\in\nbvl$ is the unique function such that
$EVg=Vh$ then we can define $g(\pm\infty)=h(\pm\infty)$.
The limit does not depend on the choice of $0\leq\lambda\leq 1$.
For more on functions of essential bounded variation, see
\cite{ambrosio} and \cite{ziemer}.

For each $n\in\N$, the operator $I^n\fn\bv\to \Ibvn$ is linear.
By the definition of $\Ibvn$ it is surjective.  Similarly if
the domain is $\nbvl$ or $\ebv$.  If $g\in\bv$ such that
$I^n[g]=0$ then by the fundamental theorem of calculus,
$\int_0^xg(t)\,dt=0$ for each $x\in\R$.  This does not imply $g=0$.
For example, $g=\chi_{\{0\}}$.  Hence, $I^n\fn\bv\to \Ibvn$
is not injective.  However, $I^n$ is a bijection when we
use $\nbvl$ or $\ebv$.

\begin{theorem}\label{theoremIbvn}
Let $n\in\N$.  Let $0\leq\lambda\leq 1$.  
(a) The sets $\Ibvn$ are equal if $\Ibv^0$ is taken
to be $\bv$, $\nbvl$ or $\ebv$.  (b) The linear operator
$I^n\fn\nbvl\to\Ibvn$ is a bijection and $\Ibvn$ is a Banach space with 
norm $\norm{h}_{\Ibvn}=
\norm{(I^n)^{-1}h}_\infty+V[(I^n)^{-1}h]$.
(c) The linear operator
$I^n\fn\ebv\to\Ibvn$ is a bijection and $\Ibvn$ is a Banach space with 
norm $\norm{h}_{\Ibvn}=
\norm{h^{(n)}}_\ebv$.
\end{theorem}

\bigskip
\noindent
{\bf Proof:} (a) Note that $\nbvl\subsetneq\bv$. 
For each element $g\in\bv$ the function
$g_\lambda\in\nbvl$ differs from $g$ on a countable set.
Hence, $I^n[g]=I^n[g_\lambda]$.  
For each $g\in\ebv$ there is exactly one function $h\in\nbvl$
such that $g=h$ almost everywhere.  Then $I^n[g]=I^n[h]$.  

(b) From (a) and Definition~\ref{defnI^n} the operator
$I^n$ is linear and surjective.
If $I^n[g]=0$ for $g\in\nbvl$ then
$\int_0^xg(t)\,dt=0$ for all $x\in\R$.  By the fundamental theorem
of calculus, $g(x)=0$ at all points of continuity of $g$.  Hence,
$g=0$ except on a countable set.  But $g$ has a left
limit and a right limit at each point.  
Suppose there is $a\in\R$ such that $g(a+)=\alpha >0$.  Then there is
$\delta>0$ such that $g(x)\geq \alpha/2$ for all $x\in(a,a+\delta)$.
This contradicts the fact that $g$ vanishes except perhaps on a
countable set. Similarly if $\alpha<0$ and similarly with the left limit.  
Hence, $g(x+)=g(x-)=0$
at all $x\in\R$.  It follows that $g=0$ on $\R$.  Hence, $I^n\fn\nbvl\to \Ibvn$
is a bijection.  If $h\in \Ibvn$ then $h=I[g]$ for a unique function
$g\in\nbvl$. The pointwise derivative $h^{(n)}(x)=g(x)$ at all points of continuity
of $g$. 
In general, we cannot recover $g(x)$ for all $x\in\R$ with
the $n$th order pointwise derivative. To compute the inverse of the 
operator $I^n$, let $S$ be the set of points in $\R$ at which $g$ is
not continuous.  Then $S$ is countable and $h^{(n)}(x)=g(x)$ for all
$x\notin S$.  Suppose $a\not\in S$. The limits
$$
\lim_{\stackrel{x\to a^+}{x\not\in S}}g(x)=g(a+)
\quad\text{ and }\quad
\lim_{\stackrel{x\to a^-}{x\not\in S}}g(x)=g(a-)
$$
both exist.  Since $g\in\nbvl$ we have $g(a)=(1-\lambda) g(a-)+\lambda g(a+)$.
This then defines $(I^n)^{-1}$.
By Lemma~\ref{lemma},
$\Ibvn$ is a Banach space with norm $\norm{h}_{\Ibvn}=
\norm{(I^n)^{-1}h}_\infty+V[(I^n)^{-1}h]$.

(c)
If $g\in\ebv$ such that
$I^n[g]=0$ then by the fundamental theorem of calculus,
$\int_0^xg(t)\,dt=0$ for each $x\in\R$.  Hence, $g=0$ almost everywhere.
But then $\langle g,\phi'\rangle=0$ for all $\phi\in\D$.
Hence, $\mu_g=0$ and $g=0$ as an element of $\ebv$.  Then
$I^n\fn\ebv\to \Ibvn$ is a bijection.  If $h\in \Ibvn$ then
the pointwise derivative $h^{(n)}(x)$ exists almost everywhere and
defines a function in $\ebv$.  By Lemma~\ref{lemma},
$\Ibvn$ is a Banach space with norm $\norm{h}_{\Ibvn}=\norm{h^{(n)}}_\ebv$.
\qed

If $h\in\Ibvn$ and $h=I^n[g]$ for $g\in\bv$ then $h=I^n[g_\lambda]$
for each $\lambda\in\R$.  Thus the normalisation on functions of
bounded variation does not affect the multiplier $h$.  If
$g\in C(\Rbar)\cap\bv$ then all normalisations $g_\lambda$ 
equal $g$ and the integral $\intinf fg$ is independent of $\lambda$.
However, if $g$ is not continuous,
different values of $\lambda$ may give different values for this
integral.  An example is given in the next section.

In the definition of $F\in\Br$ the condition
$\lim_{x\to-\infty}F(x)=0$ is imposed.  
This is arbitrary but convenient because it makes primitives 
unique.  If we merely require $\lim_{x\to-\infty}F(x)$ to exist in $\R$
then formula \eqref{intdefn2} 
must be modified by the addition of the term $(-1)^{n}F(-\infty)g(-\infty)$.
Since an element of $\arn$ is the $n$th order derivative
of a function in $\Br$, adding 
a polynomial of degree at most $n-1$
does not affect $\intinf F^{(n)}h$ or $\arn$.  The norm on $\Br$ could
then be modified to a difference formula.  Define $J_n[F](h;x)=\sum
_{i=0}^n(-1)^n\binom{n}{i}F(x+ih)$.  Then $J_n[P]=0$ if $P$ is a
polynomial of degree at most $n-1$.  The norm on $\Br$ could then
be replaced by $\sup_{x,h\in\R}|J_n[F](h;x)|$.  For example, if
$C(\Rbar)$ is used
instead of $\Bc$ then
use the norm
$\sup_{x,h\in\R}|F(x)-F(x+h)|$. 

In the definition of $h\in\Ibvn$ we have arbitrarily imposed the condition
$h^{(m)}(0)=0$ for all $0\leq m\leq n-1$.  This does not affect the
integral $\intinf fh$ if we use \eqref{intdefn2} to define the integral.  
Different lower limits of integration in
Definition~\ref{defnI^n} would change $h$ by the addition of a polynomial
of degree at most $n-1$.  Addition of such a polynomial also does
not affect the norm $\norm{\cdot}_{\Ibvn}$.
\begin{prop}\label{propmoments}
(a) Let $n\geq 2$.
Let 
$f\in\arn$ such that $f=F^{(n)}$ for $F\in\Br$. For each polynomial $P$ of degree at most
$n-2$, define $\intinf fP=(-1)^{n-1}\intinf F'P^{(n-1)}$.  Then $\intinf fP=0$.  
In particular, if $P_k(x)=x^k$ then each
of the moments $\intinf fP_k=0$ for all integers $0\leq k\leq n-2$. 
(b) Let $g\in\bv$. Let $a_1, a_2,\ldots, a_{n-1}$ be real numbers.
Define $h(x)=\lint_{x_n=a_{n-1}}^x\cdots\lint_{x_i=a_i}^{x_{i+1}}\cdots\lint
_{x_1=a_1}^{x_2}g(x_1)\,dx_1\cdots dx_{i}\cdots dx_n$.  With $f\in\arn$ and 
$F\in\Br$ such
that $f=F^{(n)}$,  define $\intinf fh
=(-1)^{n-1}\intinf F'h^{(n-1)}$.  Then $\intinf fh=\intinf fI^{n-1}[g]$.
(c) Let $n\geq 1$.  Let $h\in\Ibvn$ and 
let $P$ be a polynomial of degree at most $n-1$.  Then
$\norm{h+P}_{\Ibvn}=\norm{h}_{\Ibvn}$.
\end{prop}

\bigskip
\noindent
{\bf Proof:} (a) Let $P$
be a polynomial of degree not exceeding $n-2$.  Then
$\intinf fP=(-1)^{n-1}\intinf F'P^{(n-1)}=0$
since $P^{(n-1)}(x)=0$ for all $x\in\R$.
(b) The difference between $h$ and $I^{n-1}[g]$ is a polynomial of degree at most
$n-2$.
(c) It follows from the proof of Theorem~\ref{theoremIbvn} that
$(I^n)^{-1}P=0$.
\qed

If $\phi$ is a test function then $\phi\in\Ibvn$ for each $n\geq 0$.
For each $f\in\arn$ with primitive $F\in\Br$ the distributional
derivative formula $\langle f,\phi\rangle=\langle F^{(n)},\phi\rangle
=(-1)^n\intinf F(x)\phi^{(n)}(x)\,dx$ agrees with the definition
of the integral in \eqref{intdefn1}.
If $F,h\in C^{n-1}(\R)$ such that $F^{(n-1)},h^{(n-1)}\in AC(\Rbar)$
then the integration by parts formula is
\begin{align*}
&\intinf F^{(n)}(x)h(x)\,dx = (-1)^n\intinf F(x)h^{(n)}(x)\,dx\\
& \quad +\sum_{k=0}^{n-1}(-1)^{n-k-1}\left[
F^{(k)}(\infty)h^{(n-k-1)}(\infty)-F^{(k)}(-\infty)h^{(n-k-1)}(-\infty)
\right].
\end{align*}
When $h\in C_c^\infty(\R)$ the limits $F^{(k)}(\pm\infty)
h^{(n-k-1)}(\pm\infty)$ vanish for each $0\leq k\leq n-1$ and each
$F$.  In the case of $F\in\Dp$ these limits are ignored in the
formula for the distributional derivative $D^nF$, even though $F^{(k)}$ 
need not have any
pointwise meaning. Similarly in the definition of the integral
\eqref{intdefn2}.  We now
show that when $F$ and $h$ have pointwise derivatives as above, 
these limits $F^{(k)}(\pm\infty)
h^{(n-k-1)}(\pm\infty)$ vanish,
provided $F^{(k)}$ has a monotonicity property.  Suppose
$h\in\Ibv^{n-1}$.  Then there is $g\in\bv$ such that $h=I^{n-1}[g]$.
Hence, 
$$
\left|h^{(n-k-1)}(x)\right|  =  \left|\int_{x_k=0}^x\cdots\int_{x_1=0}^{x_2}
g(x_1)\,dx_1\cdots dx_k\right|
  \leq  \frac{\norm{g}_\infty |x|^k}{k!}
$$
so $h^{(n-k-1)}(x)=O(x^k)$ as $x\to\infty$.  This growth condition is
sharp; take $g$ to be constant.  Suppose $F\in\Bc$ such that
$F\in C^{k-1}(\R)$, $F^{(k-1)}\in AC(\Rbar)$ and there is $M>0$ such
that $F^{(k)}(x)>0$ and $F^{(k)}(x)$ is decreasing for almost all $x>M$.
Then $F$ is given by the iterated integral
$$
F(x)=\int_{x_k=-\infty}^x\cdots\int_{x_1=-\infty}^{x_2}F^{(k)}(x_1)\,dx_1\cdots dx_k.
$$
Now consider 
$I(x):=\int_{x_k=x/2}^x\cdots\int_{x_1=x/2}^{x_2}F^{(k)}(x_1)\,dx_1\cdots dx_k.$
For large enough $x$ we
have $I(x)\geq F^{(k)}(x)x^k/(2^kk!)$.
Since $\lim_{x\to\infty}I(x)=0$ it follows that $F^{(k)}(x)=o(x^{-k})$
as $x\to\infty$.  But then $\lim_{x\to\infty}F^{(k)}(x)h^{(n-k-1)}(x)=0$.
If $F^{(k)}(x)$ is increasing, instead of $I(x)$ integrate over
the interval
$[x,2x]^k\subset\R^k$.
Similarly with limits as $x\to-\infty$. 

\section{Examples and properties of the integral}\label{sectionexamples}
The space $\ac^1$ consists of the derivative of functions in $\Bc$.
Hence, it contains all functions integrable in the Lebesgue, 
Henstock--Kurzweil and wide Denjoy sense over $\R$.  For each interval
$I\subset\Rbar$ the characteristic function $\chi_I$ is of bounded
variation.  So if $f\in\ac^1$ with primitive $F\in\Bc$ then
$\int_a^bf=\intinf f\chi_I=F(b)-F(a)$.  
The same formula holds for integration over
$I=[a,b], [a,b), (a,b]$ and $(a,b)$.  Similarly, we can integrate
over all semi-infinite intervals.  If
$F$ is continuous on $\Rbar$ but has a pointwise derivative 
nowhere then $F'\in\ac^1$ and $\int_a^bF'=F(b)-F(a)$ for all
$-\infty\leq a<b\leq \infty$.  If $F\in\Bc$ is a continuous singular
function, $F'(x)=0$ for almost all $x\in\R$, then the Lebesgue
integral $\int_a^b F'(x)\,dx=0$ but $\int_a^bF'=F(b)-F(a)$.  Other
examples of integration in $\ac^1$ are given in \cite{talviladenjoy}.

The Schwartz space, $\Sc$, of rapidly decreasing test functions, consists of 
the functions $\psi\in C^\infty
(\R)$ such that for each $m,n\geq 0$ we have $x^m\psi^{(n)}(x)\to 0$
as $|x|\to\infty$.  Let $\psi\in\Sc$.  Then $\psi^{(m)}\in\Bc$ for
each $m\geq 0$.  Define $\Psi(x)=\int_{-\infty}^x\psi(t)\,dt$.  Then
$\Psi^{(m)}\in\Bc$ for each $m\geq 0$ and $\Psi'(x)=\psi(x)$ for each
$x\in\R$.
For each $1\leq n\leq m+1$, $\psi^{(m)}\in\acn$.  An example of a 
function in $\Sc$ is $\psi(x)=\exp(-x^2)$.   If we take $F_1(x)
=\exp(-x^2)$ and $F_2(x)=\int_{-\infty}^x\exp(-t^2)\,dt$ then
$F_1, F_2\in\Bc$ and if $f(x)=F_1'(x)=F_2''(x)=-2x\exp(-x^2)$
then $f\in\ac^1\cap\ac^2$. Note that $\norm{f}_{a,1}=\norm{F_1}_\infty
=1$ while $\norm{f}_{a,2}=\intinf\exp(-t^2)\,dt=\sqrt{\pi}$, so a
distribution can have different norms in different spaces $\acn$.

The space $\ar^1$ consists of the distributional derivative of
regulated functions.  Clearly, $\ac^1\subsetneq\ar^1$.  The Dirac
distribution is $\delta=H_0'$.  
Hence, $\delta\in\ar^1$.
According to \eqref{ar1defn},
for each $g\in\bv$ we have,
\begin{eqnarray*}
\intinf \delta g & = & H_0(\infty)g(\infty)-\intinf H_0(x)\,dg(x)\\
 & & \quad -[H_0(0)-H_0(0+)][g(0)-g(0+)]\\
 & = & g(\infty)-[H_0(\infty)g(\infty)-H_0(0+)g(0+)]-[0-1][g(0)-g(0+)]\\
 & = & g(0).
\end{eqnarray*}
The Henstock--Stieltjes integral $\intinf H_0(x)\,dg(x)$
can be evaluated using a tagged partition of $\Rbar$ that forces
$0$ and $\infty$ to be tags.
This agrees
with the action of $\delta$ as a tempered distribution, for which
$g$ must be in the Schwartz space ${\cal S}$.
Notice that changing the value of $H_0(0)$ does not affect the value
of $\intinf \delta g$.
When $\delta$ acts as a measure, this equation is written
$\intinf g(x)\,d\delta(x)=g(0)$ and holds for all functions
$g\fn\R\to\R$.  
Similarly, every signed Radon measure is in $\ar^1$.  

Note that
changing a function of bounded variation at one point can affect
the value of $\intinf fg$.  For example, let $F\in\Br$ and $g=a\chi
_{\{0\}}$.  Then $\intinf F'g = a[F(0+)-F(0-)]$.  See \cite{talvilaregulated}
for more examples of integration in $\ar^1$.

Using \eqref{intdefn4}, an example in $\arn$ is 
\begin{eqnarray*}
\intinf \delta^{(n-1)}I^{n-1}[H_\lambda] & = &
\intinf H_0^{(n)}I^{n-1}[H_\lambda]   =   
(-1)^{n-1}\intinf H_\lambda(x)\,dH_0(x)\\
 & = & (-1)^{n-1} H_\lambda(0)[H_0(0+)-H_0(0-)]\\
 & = &  (-1)^{n-1} \lambda.
\end{eqnarray*}
Hence, the choice of $\lambda$ affects the value of the integral.

\begin{prop}\label{propdelta}
Let $m\geq 0$.  (a) Let $n\geq 1$.  Then $\delta^{(m)}\in\arn$ if and
only if $n=m+1$.  (b) For no $n\geq 1$ is $\delta^{(m)}\in\acn$.
(c) $\norm{\delta^{(m)}}_{a,m+1}=1$.
(d) Let $\mu$ be a finite signed Borel measure.  Then $D^m\mu\in\ar^{m+1}$.
\end{prop}

\bigskip
\noindent
{\bf Proof:} (a) 
We have $\delta^{(m)}=D^{m+1}H_0$, so $\delta^{(m)}\in\ar^{m+1}$.
Suppose $\delta^{(m)}= F^{(n)}$
for some $n\geq m+2$.  Then $F^{(n-m)}=\delta+P$ where $P$ is a polynomial
of degree at most $m-1$ and the equality is in $\ar^1$.  If $m=0$ then
$P=0$.  Integrating
$n-m$ times over the interval $[0,x]$ gives 
$F(x)=x^{n-m-1}\chi_{(0,\infty]}(x)/(n-m-1)!+Q(x)$ where $Q$ is a
polynomial of degree at most $n-1$.  But then $F\not\in\Br$.  If $n\leq m$ then
$F=\delta^{(m-n)}+P$ where $P$ is a polynomial of degree at most $n-1$.
Comparing supports shows this impossible for all $F\in\Br$.
(b) Part (a) includes the proof. 
(c)  Notice that $\norm{\delta^{(m)}}_{a,m+1}=\norm{H_0}_\infty
=1$ for each $m\geq 0$.
(d) Define $F(x)=\int_{(-\infty,x)}d\mu$.  Then $F\in\nbvzero\subset\Br$ and
$F'=\mu$. \qed

Observe that
\eqref{intdefn1} gives
$\intinf \delta^{(m)}I^{m}[g]=(-1)^mg(0)$, for
each $g\in\nbvl$ and $m\geq 0$.
Let $K(x)=x\chi_{[0,\infty]}(x)$.  Then
$\delta^{(m)}=D^{m+2}K$ for each $m\geq 0$.  Although $K$ is continuous,
it is not in $\Bc$.  However, if we let $F_1(x)=0$ for $x\leq 0$,
$F_1(x)=x$ for $0\leq x\leq 1$, $F_1(x)=1$ for $x\geq 1$ then
$F_1\in\Bc$. Let $F_2(x)=H_0(x)-H_0(x-1)$.  Then $F_2\in\Br\setminus
\Bc$.  Let $f=F_1''\in\ac^2$ then $f=F_2'\in\ar^1$ and $f=\delta
-\tau_1\delta$. So a linear combination of elements from $\ar^1$
is in $\ac^2$.

\begin{prop}
(a) For each $1\leq m<n$, $\acn$ is not a subset of $\ac^m$ and
$\ac^m$ is not a subset of $\acn$.
(b) $\ac^m\subset\arn$ if and only if $m=n$.
(c) For each $1\leq m<n$, $\arn$ is not a subset of $\ar^m$ and
$\ar^m$ is not a subset of $\arn$.
(d) For no $m,n\in\N$ is $\ar^m\subset\acn$.
(e) For each $m,n\in\N$, $\ac^m\cap\acn\not=\varnothing$.  Hence,
$\ar^m\cap\arn\not=\varnothing$.
\end{prop}

\bigskip
\noindent
{\bf Proof:} (a)  There is an increasing function $F\in\Bc\cap C^\infty(\R)$ 
such
that $F=0$ on $(-\infty,0]$ and $F=1$ on $[1,\infty)$.  Let $f=F^{(m)}\in\ac^m$.
Suppose $f\in\acn$.
The only function $G\in\Bc$ that satisfies
$G^{(n)}=f$ is given by the iterated improper Riemann integrals
\begin{eqnarray*}
G(x) & = & \lint_{x_n=-\infty}^x\cdots\lint_{x_i=-\infty}^{x_{i+1}}\cdots\lint
_{x_1=-\infty}^{x_2} F^{(m)}(x_1)\,dx_1\cdots dx_n\\ 
 & = & \lint_{x_{n}=-\infty}^x\cdots\lint_{x_i=-\infty}^{x_{i+1}}\cdots\lint
_{x_{m+1}=-\infty}^{x_{m+2}} F(x_{m+1})\,dx_{m+1}\cdots dx_{n}\\ 
 & \geq & \lint_{x_{n}=1}^x\cdots\lint_{x_i=1}^{x_{i+1}}\cdots\lint
_{x_{m+1}=1}^{x_{m+2}} \,dx_{m+1}\cdots dx_{n}\quad \text{ if } x\geq 1\\ 
 & = & \frac{(x-1)^{n-m}}{(n-m)!}.
\end{eqnarray*}
Hence, $G\not\in\Bc$ so $f\not\in\acn$.

Let $F(x)=0$ for $x\leq 0$, $F(x)=x$ for $0\leq x\leq 1$ and
$F(x)=1$ for $x\geq 1$.  Then $F\in\Bc$.  Define $f=F^{(n)}\in\acn$.
The Heaviside step function is $H_0=\chi_{(0,\infty]}$.  The Dirac
distribution is $\delta=H_0'$.  If $a\in\R$ we write $\delta_a=
\tau_a\delta$ for the Dirac distribution supported at $a$.
For $n\geq 2$
we have
$f=F^{(n)}=\delta^{(n-2)}-\delta_1^{(n-2)}$.  Suppose $f\in\ac^m$
is
given by $f=G^{(m)}$ for $G\in\Bc$.
If $n\geq m+2$ then
$G=\delta^{(n-m-2)}-\delta_1^{(n-m-2)}+P$ where
$P$ is a polynomial of degree at most $m-1$.
If $n=m+1$ then $G(x)=H_0(x)-H_0(x-1)$.
It follows that
$G\not\in\Bc$.  Hence, $f\not\in\ac^m$.

(b) Since $\Bc\subset\Br$ we have $\acn\subset\arn$.  For the
other part of the proof use examples as in part (a).
Replace the second example in (a) by $F(x)=H_0(x)-H_0(x-1)$.

(c) By Proposition~\ref{propdelta}, $\delta^{(m-1)}\in\ar^n$
if and only if $n=m$.

(d)  By Proposition~\ref{propdelta}, $\delta^{(m-1)}\in\ar^m$
but is not in any of the $\acn$ spaces.

(e) See the example in the second paragraph of this
section.\qed

Let $F(x)=H_0(x-1)(x-1)^\alpha e^{-(x-1)}$ where $\alpha>0$.
Then $F\in\Bc$.  Define $f\in\acn$ by $f=F^{(n)}$.  
For all $x\not=1$ the pointwise derivative gives
$f(x)\sim\alpha(\alpha-1)\cdots(\alpha-n+1)(x-1)^{\alpha-n}$
as $x\to 1+$.  Then $f$ is singular at $1$ such that $f\not\in
L^1_{loc}$ if $\alpha\leq n-1$ and yet $\intinf fh$
exists for each $h\in\Ibvn$.

In $\acn$ there is a version of the second mean value theorem for 
integrals.
\begin{theorem}
Let $F\in\Bc$.  Let $h\in\Ibv^{n-1}$ such that $h^{(n-1)}$ is a
monotonic function.  Then
$\intinf F^{(n)}h=(-1)^{n-1}[h^{(n-1)}(-\infty)\int_{-\infty}^\xi F'
+h^{(n-1)}(\infty)\int_\xi^\infty F']$ for some $\xi\in\Rbar$.
\end{theorem}

\bigskip
\noindent
{\bf Proof:} Integrate and use the mean value theorem for
Riemann--Stieltjes integrals \cite[\S7.10]{mcleod}:
\begin{eqnarray*}
\intinf F^{(n)}h & = & (-1)^{n-1}\left[F(\infty)h^{(n-1)}(\infty)-\intinf 
F(x)\,dh^{(n-1)}(x)\right]\\
 & = & (-1)^{n-1}\left[F(\infty)h^{(n-1)}(\infty)-F(\xi)\intinf dh^{(n-1)}(x)
\right]\\
 & = &  (-1)^{n-1}\left\{F(\infty)h^{(n-1)}(\infty)-F(\xi)[h^{(n-1)}(\infty)-
h^{(n-1)}(-\infty)]\right\}\\
 & = & (-1)^{n-1}\left\{h^{(n-1)}(-\infty)F(\xi) + h^{(n-1)}(\infty)[F(\infty)-F(\xi)]
\right\}.\qed
\end{eqnarray*}
This proof is adapted from a similar theorem for the wide Denjoy integral
in \cite{celidze}, where a proof of the Bonnet
form of the second mean value theorem can also be found.

If $\psi\in C^\infty(\R)$ is a bijection such that $\psi'>0$ on $\R$
then for any distribution $T\in\Dp$ the composition $T\circ\psi$
is defined by $\langle T\circ\psi,\phi\rangle=\langle T,\tfrac{\phi
\circ \psi^{-1}}{\psi'\circ\psi^{-1}}\rangle$ for all $\phi\in\D$.
In \cite{talviladenjoy} a change of variables formula was proved in
$\ac^1$ when
$\psi\in C(\R)$, i.e., no monotonicity or pointwise differentiability
is assumed.  In \cite{talvilaregulated} a change of variables formula 
was proved in
$\ar^1$ when
$\psi$ was piecewise monotonic.  
For $\arn$ we have the simple case of composition
with a linear function.
\begin{theorem}\label{theoremchangeofvariables}
Let $\psi(x)=ax+b$ for $a,b\in\R$, $a\not=0$.  Let $F\in\Br$.  Then
$(F\circ \psi)^{(n)}=a^n(F^{(n)}\circ\psi)$.  Let $h\in\Ibv^{n-1}$.
Then $\intinf F^{(n)}h=|a|\intinf (F^{(n)}\circ\psi)(h\circ\psi)$.
\end{theorem}

\bigskip
\noindent
{\bf Proof:} Let $\phi\in\D$.  Then
\begin{eqnarray*}
\left\langle F^{(n)}\circ\psi,\phi\right\rangle & = & 
{\rm
sgn}(a)\left\langle F^{(n)},\frac{
\phi
\circ \psi^{-1}}{\psi'\circ\psi^{-1}}\right\rangle\\
 & = & 
\frac{{\rm sgn}(a)(-1)^n}{a^{n+1}}\intinf F(y)\phi^{(n)}\circ\psi^{-1}(y)\,dy\\
 & = &  \frac{(-1)^n}{a^{n}}\intinf F(ax+b)\phi^{(n)}(x)\,dx\\
 & = & a^{-n}\langle(F\circ\psi)^{(n)},\phi\rangle.
\end{eqnarray*}
This shows that
$(F\circ\psi)^{(n)}=a^n(F^{(n)}\circ\psi)$.  Note that $F\circ\psi\in\Br$.
Suppose $a>0$.   Then
\begin{align}
&a^n\intinf(F^{(n)}\circ\psi)\,h\circ\psi  =  \intinf (F\circ\psi)^{(n)}h\circ
\psi\notag\\
&=(-1)^{n-1}a^n\intinf h^{(n-1)}(ax+b)\,dF(ax+b)\notag\\
&=(-1)^{n-1}a^{n-1}\intinf h^{(n-1)}(y)\,dF(y)\label{stieltjeschange3}\\
&=a^{n-1}\intinf F^{(n)}h.\notag
\end{align}
If $a<0$ then there is a sign change in \eqref{stieltjeschange3}
upon change of variables.\qed

Define $r_x(y)=x-y$.  Then for $F\in\Br$ and $h\in\Ibv^{n-1}$,
Theorem~\ref{theoremchangeofvariables} shows the equality of the
two convolution
integrals $\intinf (F^{(n)}\circ r_x)h=\intinf F^{(n)}(h\circ r_x)$.

If $f\in\acn$ for $n\geq 2$ then in general $\int_a^b f$ does not
exist.  However, if $h\in\Ibv^{n-1}$ has compact support and is
in $C^{n-1}(\R)$ then $fh$ can be integrated over a subinterval.
For example, let $a<b$.  Define $h(x)=(x-a)^p(x-b)^q\chi_{[a,b]}(x)$
for $p,q\geq n-1$.  There is a polynomial, $P$, of degree at most
$n-2$ such that $h+P\in\Ibv^{n-1}\cap C^{n-2}(\R)$.  
It follows that $h+P=I^{n-1}[g]$
where 
$$
g(x)=\chi_{(a,b)}(x)\lsum_{i=0}^{n-1}\!\binom{n-1}{i}(p-i+1)_i(q-n+i+2)_{n-i-1}
(x-a)^{p-i}(x-b)^{q-n+i+1}.
$$
The Pochhammer symbol is $(z)_m=z(z+1)\cdots(z+m-1)$ for $m\in\N$ with
$(z)_0=1$.  The formula for $g$ comes from the Leibniz rule for differentiating
a product.  Note that $g\in\bv$.  The value of $g$ at $a$ and $b$ is
irrelevant.  Let $f=F^{(n)}$ for
$F\in\Bc$.
Then since $g(a-)=g(b+)=0$, Proposition~\ref{propmoments} gives
\begin{eqnarray*}
\intinf fh & = & (-1)^{n-1}\intinf F' h^{(n-1)}=(-1)^{n-1}\intinf F'g\\
 & = & (-1)^{n}\left[F(a)g(a)-F(b)g(b)+\int_a^b F(x)\,dg(x)\right].
\end{eqnarray*}
This defines
$\int_a^b F^{(n)}g=
(-1)^{n}[F(a)g(a)-F(b)g(b)]+(-1)^{n}\int_a^b F(x)\,dg(x)$.
If $p,q>n-1$ then $h\in C^{n-1}(\R)$, $g$ is continuous, $g(a)=g(b)=0$ and 
$\intinf fh=(-1)^n\int_a^bF(x)
\,dg(x)$.  
If $p,q>n$ then $g\in C^1(\R)$ and from Proposition~\ref{propLebesgueRiemann},
$\intinf fh=(-1)^n\int_a^bF(x)g'(x)\,dx$.  If $f\in\arn$ then we must
adopt a normalisation on $g$.  If $g\in\nbvl$ then
\begin{eqnarray*}
\intinf fh & = & (-1)^{n}\left\{F(a)[g_\lambda(a)-g_\lambda(a-)]
+F(b)[g_\lambda(b+)-g_\lambda(b)]\right.\\
 & & \quad+\int_a^b F(x)\,dg_\lambda(x)
+[F(a)-F(a+)][g_\lambda(a)-g_\lambda(a+)]\\
 & & \quad\left.+[F(b)-F(b+)][g_\lambda(b)-g_\lambda(b+)]\right\}.
\end{eqnarray*}

There are functions in $\Ibvn$ that play the role
of characteristic functions of intervals.  These lead to a version
of the fundamental theorem of calculus that is built in to the 
definition of the integral.  This also gives an explicit formula
for the inverse of the $n$th derivative operator.

\begin{theorem}\label{theoremftc}
Let $x\in\R$.  Let $\lambda=1$.
Define $h(x,t)=(x-t)^{n-1}H_0(x-t)/(n-1)!$.  (a) $h(0,\cdot)=I^{n-1}[
(-1)^{n-1}\chi_{[-\infty,0)}]\in\Ibv^{n-1}$.  (b) Let $F\in\Br$.  Then
$\intinf F^{(n)}h(x,\cdot)
=F(x)$.  (c) Let $f\in\arn$ and define $G(x)=\intinf fh(x,\cdot)$.  Then
$G^{(n)}=f$. (d) The operator $\Phi\fn\arn\to\Br$ defined by
$\Phi[f](x)=\intinf fh(x,\cdot)$ is a linear isometry and is the inverse
of $D^n\fn\Br\to\arn$ given by $D^n[F]=F^{(n)}$.
Similarly with $\Phi\fn\acn\to\Bc$.
\end{theorem}

\bigskip
\noindent
{\bf Proof:} (a) Let $g(x,t)=(-1)^{n-1}\chi_{[-\infty,x)}(t)$.  Note
that $g(x,\cdot)\in \nbvone$, $h(x,\cdot)\in C^{n-2}(\R)$ and that 
$\partial^{n-1}h(x,t)/\partial t^{n-1}
=g(x,t)$
for all $t\not=x$.  Therefore, $h(0,\cdot)=I^{n-1}[g(0,\cdot)]\in \Ibv^{n-1}$.
(b) From \eqref{intdefn1} and Proposition~\ref{propmoments},
$\intinf F^{(n)}h(x,\cdot)=(-1)^{n-1}\intinf F' (-1)^{n-1}\chi_{[-\infty,x)}
=\int_{[-\infty,x)} F'=F(x-)=F(x)$.  (c) This follows from (b).\qed

Comparing the result of Proposition~\ref{propIbvnint} and
Theorem~\ref{theoremftc}(b), it is clear how to define the iterated
integral of the $n$th derivative of functions in $\Bc$.
\begin{defn}\label{defniterated}
Let $F\in\Bc$.
Define $\lint_{x_n=-\infty}^x\cdots\lint_{x_i=-\infty}^{x_{i+1}}\cdots\lint
_{x_1=\-\infty}^{x_2}F^{(n)}=F(x)$ for each $x\in\R$.
\end{defn}
The definition can be justified as follows.  The set of functions
$\Bc^\infty$ is defined to be those functions $\psi\in C^\infty(\R)$ 
for which there are real numbers $a<b$ and $c$ such that $\psi=0$ on 
$(-\infty,a]$ and $\psi=c$ on $[b,\infty)$.  It is clear that
$\Bc^\infty$ is dense in $\Bc$.  If $(\psi_k)$ is a Cauchy sequence
in $\Bc^\infty$ then
\begin{align*}
&\left|\,\lint_{x_n=-\infty}^x\!\!\cdots\!\!\lint
_{x_1=-\infty}^{x_2}\!\psi_k^{(n)}(x_1)\,dx_1\cdots  dx_n
-\lint_{x_n=-\infty}^x\!\!\cdots\!\!\lint
_{x_1=-\infty}^{x_2}\!\psi_l^{(n)}(x_1)\,dx_1\cdots dx_n\right|\\
&=|\psi_k(x) -\psi_l(x)|\leq\norm{\psi_k-\psi_l}_\infty.
\end{align*}
Hence, the sequence of functions $\lint_{x_n=-\infty}^x\!\!\cdots\!\!\lint
_{x_1=-\infty}^{x_2}\!\psi_k^{(n)}(x_1)\,dx_1\cdots  dx_n$ has
a limit it $\Bc$.  We can define the limit to be $F(x)$ in the case
when $F\in\Bc$ and $\lim_{k\to\infty}\norm{F-\psi_k}_\infty=0$.  It is
easy to see that the value of the limit is independent of the choice
of sequence $(\psi_k)$.

It is a classical result that 
the initial value problem; given $f\in C(\R)$ such that $\int_{-\infty}^0
|f(t)||t|^{n-1}\,dt<\infty$, find $F\in C^n(\R)$ such
$F^{(n)}(x)=f(x)$ for all $x\in\R$, with initial condition
$\lim_{x\to-\infty}F^{(k)}(x)=0$ for each $0\leq k\leq n-1$,
has the unique solution $F(x)=[1/(n-1)!]\int_{-\infty}^x f(t)
(x-t)^{n-1}\,dt$.  By the Fubini--Tonelli theorem the solution
can also be written $F(x)=\int_{-\infty}^x\cdots\int_{-\infty}^{x_2}
f(x_1)\,dx_1\cdots dx_n$.  Hence, an alternative approach to the integral
is to use Theorem~\ref{theoremftc}(b) and Definition~\ref{defniterated},
rather 
than Definition~\ref{defnint}.

\section{H\"older inequality and dual space}\label{sectionholder}
One of the many useful properties of functions in an $L^p$ space
is the H\"older inequality.  Distributions in $\arn$
also satisfy a type of H\"older inequality.  For each
$f\in\arn$ and each $h\in\Ibv^{n-1}$ the integral Definition~\ref{defnint}
provides a type of product.

\begin{theorem}[H\"older inequality]
Let $f\in\arn$ with primitive $F\in\Br$.  Let $h\in\Ibv^{n-1}$ such that
$h=I^{n-1}[g]$ for $g\in\bv$.  Then 
$|\intinf fh|\leq \norm{F}_\infty\norm{g}_\bv=\norm{f}_{a,n}\norm{h}_{
\Ibv^{n-1}}$.
\end{theorem}
The proof follows from Definition~\ref{defnint}.  The case $n=1$ 
was considered in \cite{talviladenjoy} for $\ac^1$ and in
\cite{talvilaregulated} for $\ar^1$ where various other forms
of this inequality can be found.  The estimates do not depend on 
the choice of $0\leq\lambda\leq 1$ since 
$\norm{g_\lambda}_\bv\leq\norm{g}_\bv$
for all $g\in\bv$.  For $g\in\ebv$ and $f\in\acn$ use $\norm{g}_\ebv$.

An application of the H\"older inequality is the following convergence
theorem.
\begin{theorem}\label{theoremconvergence}
Fix $n\in\N$.  Fix $0\leq\lambda\leq 1$.
Let $f\in\arn$ and
for each $k\in\N$
let $f_k\in\arn$ such that $\norm{f_k-f}_{a,n}\to0$. 
Let $h,h_k\in\Ibv^{n-1}$ such that $\norm{h_k-h}_{\Ibv^{n-1}}\to0$.
Then $\lim_{k\to\infty}\intinf f_kh_k=\intinf fh$.
\end{theorem}

\bigskip
\noindent
{\bf Proof:} 
The H\"older inequality gives
\begin{eqnarray*}
\left|\intinf f_kh_k -\intinf fh\right| & \leq & 
\left|\intinf (f_k -f)h_k\right|+ \left|\intinf f(h_k-h)\right|\\
 & \leq & \norm{f_k-f}_{a,n}\norm{h_k}_{\Ibv^{n-1}}+\norm{f}_{a,n}\norm
{h_k-h}_{\Ibv^{n-1}}.
\end{eqnarray*}
Since $\norm{h_k}_{\Ibv^{n-1}}$ is bounded, the result follows.\qed

The H\"older inequality shows that for each $f\in\acn$ the integral is a 
continuous linear
functional on $\Ibv^{n-1}$ and that for each $h\in\Ibv^{n-1}$ the integral is a 
continuous linear
functional on $\acn$.  There is also an equivalent norm in terms of
these functionals.  Similarly for $\arn$.

\begin{theorem}
(a) Let $f\in\acn$.  Define $\Phi_f\fn\Ibv^{n-1}\to\R$ by
$\Phi_f(h)=\intinf fh$.  Then $\Phi_f$ is a continuous linear functional.
(b) Let $h\in\Ibv^{n-1}$.  Define $\Psi_h\fn\Ibv^{n-1}\to\R$ by
$\Psi_h(f)=\intinf fh$.  Then $\Psi_h$ is a continuous linear functional.
(c) Let $f\in\acn$. Define $\norm{f}'_{a,n}=\sup_h\intinf fh$ where the supremum is taken
over all $h\in\Ibv^{n-1}$ with $h=I^{n-1}[g]$ for $g\in\bv$ such that
$\norm{g}_\infty\leq 1$ and $Vg\leq 1$.  Then $\norm{\cdot}_{a,n}$
and $\norm{\cdot}'_{a,n}$ are equivalent norms on $\acn$.
(d) Fix $0\leq\lambda\leq 1$.  Results analogous to (a), (b) and (c)
hold for $\arn$.
\end{theorem}

\bigskip
\noindent
{\bf Proof:} Linearity in each argument follows from linearity of
the derivatives defining $f$ and the integrals defining $h$.
If $f\in\acn$ and $(h_k)\subset\Ibv^{n-1}$ such that $\norm{h_k}_{\Ibv^{n-1}}\to 0$
as $k\to\infty$  then by Theorem~\ref{theoremconvergence}, $\Phi_f(h_k)\to 0$.
If $h\in\Ibv^{n-1}$ and $(f_k)\subset\acn$ such that $\norm{f_k}_{a,n}\to 0$
as $k\to\infty$ then Theorem~\ref{theoremconvergence}, $\Psi_h(f_k)\to 0$.
Part (c) follows from Theorem~$29$ in \cite{talviladenjoy}, which
proves equivalence on $\ac^1$. (d) See \cite[Theorem~15]{talvilaregulated}
for equivalent norms in $\arn$.\qed

It is a classical result that if $[a,b]$ is
a compact interval then $C([a,b])^\ast=\nbvl$ and if $L$ is an element of
the dual space then there is a function $g\in\nbvl$ such that 
$L(F)=\int_a^b F(x)\,dg(x)$ for all $F\in C([a,b])$.
By the compactification, the same holds for $C(\Rbar)$.
The choice of $0\leq\lambda\leq 1$ is immaterial.
The distinction between
$\bv$ and $\nbvl$ is sometimes ignored (including in \cite{talviladenjoy}
in the paragraphs preceding Theorem~8).  The distinction is important since if
$F$ is continuous and $g=\chi_{\{0\}}$ then $\intinf F(x)\,dg(x)=0$.
This
function $g$ is not $0$ as an element of $\bv$ but its normalised
version in $\nbvl$ is $0$.  It was shown in \cite{talviladenjoy} that
the dual space of $\ac^1$ is $\ebv$.  It then follows from Lemma~\ref{lemma}
that the dual space
of $\acn$ is also $\ebv$.  For each $0\leq\lambda\leq 1$  we can choose
to represent each element
of $\ebv$ by a unique function in $\nbvl$.  Hence, the dual space of
$\acn$ is $\nbvl$.  Meanwhile, the dual space of $\Br$ is $\bv$.  See
\cite{talvilaregulated} for a discussion of this point.
\begin{theorem}
For each $n\in\N$ the dual space of $\acn$ 
is isometrically isomorphic
to $\ebv$ and $\nbvl$.  The dual space of $\arn$ is isometrically
isomorphic to $\bv$.
\end{theorem}

By Theorem~\ref{theoremIbvn} and Lemma~\ref{lemma}, the spaces
$\nbvl$ and $\Ibvn$ are isometrically isomorphic.  Hence, $\acn\subset
\nbvl^\ast$.
Spaces of distributions are often defined as the dual of some topological
vector space.  No explicit description of the dual of $\bv$ or $\nbvl$ seems
to be known.  We are thus reluctant to define our space of integrable
distributions as $\nbvl^\ast$ and instead have chosen the concrete 
description in terms of derivatives of functions in $\Bc$ and $\Br$.

\section{Banach lattice}\label{sectionlattice}
The usual pointwise ordering
makes $\Br$ into a Banach lattice.  Each of the spaces
$\arn$ inherits this Banach lattice structure.  We will point out a
few of the most basic lattice properties of $\arn$ but leave a detailed study
for later.  A reference for this section is \cite{aliprantisborder}.
To keep the paper reasonably self contained we prove all results in this
section ab initio, although some of them follow from more general lattice
theorems.

If $\preceq$ is a binary
operation on set $S$ then it is a {\it partial order} if for all
$x,y,z\in S$ it is {\it reflexive} ($x\preceq x$), {\it antisymmetric}
($x\preceq y$ and $y\preceq x$ imply $x=y$) and {\it transitive} ($x\preceq y$
and $y\preceq z$ imply $x\preceq z$). 
If
$S$ is a Banach space with norm $\norm{\cdot}_S$ and $\preceq$ is a partial
order on $S$  then $S$ is a {\it Banach lattice} if for all $x,y,z\in S$
\begin{enumerate}
\item
$x\vee y$ and $x\wedge y$ are in $S$.  The {\it join} is 
$x\vee y=
\sup\{x,y\}=w$ such that $x\preceq w$, $y\preceq w$ and if $x\preceq\tilde{w}$
and $y\preceq\tilde{w}$ then $w\preceq\tilde{w}$.
The {\it meet} is 
$x\wedge y=
\inf\{x,y\}=w$ such that $w\preceq x$, $w\preceq y$ and if $\tilde{w}\preceq x$
and $\tilde{w}\preceq y$ then $\tilde{w}\preceq w$.
\item
$x\preceq y$ implies $x+z\preceq y+z$.
\item
$x\preceq y$ implies $kx\preceq ky$ for all $k\in\R$ with $k\geq 0$.
\item
$|x|\preceq |y|$ implies $\norm{x}_S\leq \norm{y}_S$.
\end{enumerate}
If $x\preceq y$ we write $y\succeq x$.
We also define
$|x|=x\vee (-x)$, $x^+=x\vee 0$ and $x^-=(-x)\vee 0$. 
Then $x=x^+-x^-$ and $|x|=x^++x^-$.

The usual pointwise ordering, $F_1\leq F_2$ if and only if $F_1(x)\leq F_2(x)$
for all $x\in\R$,  is a partial order on $\Br$.
Since $\Br$ is closed under the operations
$(F_1\vee F_2)(x)=\sup(F_1,F_2)(x)=\max(F_1(x),F_2(x))$ and
$(F_1\wedge F_2)(x)=\inf(F_1,F_2)(x)=\min(F_1(x),F_2(x))$,
it is then a vector lattice (or Riesz space).
The inequality $\norm{F_1F_2}_\infty\leq\norm{F_1}_\infty\norm{F_2}_\infty$
shows
$\Br$ is also a Banach lattice.  See \cite{talvilaregulated}.  Clearly,
$\Bc$ is a sublattice.  Notice that the ordering in $\Br$ depends on
our choice of using left continuous primitives.

A partial ordering in $\arn$ is inherited from $\Br$.  If $f_1, f_2\in
\acn$ with respective primitives $F_1, F_2\in\Br$ then  $f_1\preceq f_2$
if and only if $F_1\leq F_2$ in $\Br$.  By Lemma~\ref{lemma},
$\arn$ is a Banach lattice and $\acn$ is a sublattice. 

An element $e\geq 0$ such that for each $x\in S$ there is $\lambda>0$
such that $|x|\leq \lambda e$ is an {\it order unit} for lattice 
$S$.  In the
theorem below we show $\Br$ and hence $\arn$ does not have an
order unit.

We have absolute
integrability: if $f\in\arn$ so  is $|f|$.  The $n$th derivative operator
$D^n$ commutes with $\vee$ and $\wedge$ and hence with $|\cdot|$.

\begin{theorem}[Banach lattice]\label{theoremlattice}
{\rm (a)} $\Br$ is a Banach lattice and $\Bc$ is a Banach sublattice.
{\rm (b)}
For $f_1,f_2\in\arn$ with respective primitives $F_1,F_2\in\Br$, 
define $f_1\preceq f_2$ if $F_1\leq F_2$ in $\Br$.
Then $\arn$ is a Banach lattice isomorphic to $\Br$.
{\rm (c)}
$\Br$ and $\arn$ do not have an order unit.
{\rm (d)}  Let $F_1,F_2\in\Br$.  Then $D^n(F_1\vee F_2)=
(F_1)^{(n)}\vee F_2^{(n)}$,
$D^n(F_1\wedge F_2)=F_1^{(n)}\wedge F_2^{(n)}$, $|F^{(n)}|=D^n|F|$,
$D^n(F^+)=(D^nF)^+$, and $D^n(F^-)=(D^nF)^-$.
{\rm (e)} If $f\in\arn$ with primitive $F\in\Br$ then $|f|\in\arn$ with 
primitive $|F|\in\Br$. Let $h=I^{n-1}[(-1)^{n-1}\chi_{(a,b)}]$ for
$(a,b)\subset\R$.  Then  
$|\intinf fh|\geq|\intinf|f|h\,|$ for all $f\in\acn$.  
Now let $h(x,t)=(x-t)^{n-1}H_0(x-t)/(n-1)!$.
Then
$|\intinf f h(x,\cdot)|
=\intinf|f|h(x,\cdot) =|F(x)|$ for all $f\in\acn$.  This formula also
holds for $f\in\arn$ if $\lambda=1$.  And, $\norm{|f|}_{a,n}=\norm{f}_{a,n}$, 
$\norm{f^{\pm}}_{a,n}\leq \norm{f}_{a,n}$.  
{\rm (f)} If $f\in\arn$ then $f^{\pm}\in\arn$ with respective
primitives $F^{\pm}\in\Br$.  {\it Jordan decomposition}: $f=f^+ - f^-$.
And, $\intinf fh=\intinf f^+h-\intinf f^-h$ for every $h\in\Ibv^{n-1}$.
{\rm (g)} $\arn$ is {\it distributive}:
$f\wedge(g\vee h)=(f\wedge g)\vee(f\wedge h)$
and  $f\vee(g\wedge h)=(f\vee g)\wedge(f\vee h)$ for all $f,g,h\in\arn$. 
{\rm (h)} $\arn$ is {\it modular}: For all $f,g\in\arn$, if
$f\preceq g$ then $f\vee(g\wedge h)=g\wedge(f\vee h)$ for all $h\in\arn$.
{\rm (i)} Let $F_1$ and $F_2$ be continuous functions on $\Rbar$. Then 
\begin{eqnarray}
F_1^{(n)}\preceq F_2^{(n)} & \Longleftrightarrow & F_1(x)-F_1(-\infty)\leq 
F_2(x)-F_2(-\infty)
\quad\forall x\in\R.\label{lattice1}
\end{eqnarray}
Let $F_1$ and $F_2$ be regulated functions on $\Rbar$.
Then 
\begin{eqnarray}
F_1^{(n)}\preceq F_2^{(n)} & \Longleftrightarrow & F(x-)-F(-\infty)\leq G(x-)-G(-\infty)
\quad\forall x\in\R\label{lattice2}\\
 &  \Longleftrightarrow & F(x+)-F(-\infty)\leq G(x+)-G(-\infty)
\quad\forall x\in\R.\label{lattice3}
\end{eqnarray}

\end{theorem}

\bigskip
\noindent
{\bf Proof:} (a) It is clear that $\Br$ is closed under supremum
and infimum.  See \cite{talvilaregulated}.  Hence, it is a Banach 
sublattice of the bounded functions on $\Rbar$.
(b) This follows from Lemma~\ref{lemma}.
(c) Suppose $e\in\Br$ is an order unit.  Then $F$ defined by
$F(x)=\sqrt{e(x)}$ is in
$\Br$.  And, $\lambda\geq\lim_{x\to-\infty}|F(x)|/e(x)=\lim_{x\to-\infty}
1/\sqrt{e(x)}=\infty$.  Hence, $\Br$ has no order unit.
This shows $\arn$ has no order unit.
(d) Suppose $f\in\arn$ with primitive $F\in\Br$ such that
$D^nF_1\vee D^nF_2=f$.  Then $D^nF_1\preceq f$, $D^nF_2\preceq f$
and if $D^nF_1\preceq\ftilde$, $D^nF_2\preceq\ftilde$ for some
$\ftilde\in\arn$ then $f\preceq\ftilde$.  These statements are
equivalent to $F_1\leq F$, $F_2\leq F$ and if $F_1\leq \Ftilde$,
$F_2\leq \Ftilde$ then $F\leq \Ftilde$, where $\Ftilde\in\Br$ is
the primitive of $\ftilde$.  Therefore,
$F=F_1\vee F_2$ so $D^n(F_1 \vee F_2)=D^nF=f=D^nF_1\vee D^nF_2$.
The other parts are similar.
(e) For $f\in\acn$, note that 
$|\intinf f h|=|\intinf F'h^{(n-1)}|=|F(b)-F(a)|$ and
$\intinf|f|h=|F(b)|-|F(a)|$.
The other parts of
(e) and (f) follow from (d) and the definitions, together with
Proposition~\ref{propmoments} and
Theorem~\ref{theoremftc}.
(g) The real-valued functions on any set form a distributed lattice
due to inheritance from $\leq$ in $\R$.  Therefore, $\Br$ is a distributed
lattice and so is $\arn$. See \cite[p.~484]{maclanebirkhoff} for an
elementary
proof and for another property of distributed lattices.  (h) Modularity
is also inherited from $\leq$ in $\R$ via $\Br$.
(i) We have $F_1^{(n)},F_2^{(n)}\in\acn$ with respective 
primitives $\Phi_1,\Phi_2\in
\Bc$ given by $\Phi_1(x)=F_1(x)-F_1(-\infty)$ and
$\Phi_2(x)=F_2(x)-F_2(-\infty)$.  The definition of order then gives 
\eqref{lattice1}.
The relations $F(x\pm)=\lim_{y\to x^\pm}F(y)=
\lim_{y\to x^\pm}F(y-)$ give \eqref{lattice2} and \eqref{lattice3}.
\qed

Let $f_1,f_2\in\arn$ with respective primitives $F_1,F_2\in\Br$.
Note that if $F_1\leq F_2$ in $\Br$ then we can differentiate both
sides of this inequality with $D^n$ to get  $f_1\preceq f_2$ in $\arn$.  And,
if $f_1\preceq f_2$ in $\arn$ we can integrate both sides against 
$h(x,\cdot)$ to get $F_1\leq F_2$ in $\Br$.  See Theorem~\ref{theoremftc}. 
This also shows the derivative $D^n$ is a positive operator on 
$\Br$ and its inverse is a positive operator on $\arn$.

Define $H_{-\infty}\fn\Rbar\to\R$ by $H_{-\infty}=\chi_{(-\infty,\infty]}$.
Then $H_{-\infty}$ behaves like an order unit for $\Br$.  However, as a
distribution it is equal to the constant distribution $1$.  Hence, all
of its distributional derivatives are $0$.  To include an order unit
we have to use a more general type of differentiation with respect to
test functions that are not necessarily $0$ at $-\infty$.  Two possibilities
are $\bv$ or $C^\infty(\R)\cap C(\Rbar)$.
We can take as a space of
primitives $\B$,
which consists of the functions $F\fn\Rbar\to\R$ that are regulated
on $\Rbar$
such that $F(-\infty)=0$, $F(x)=F(x-)$ for all $x\in\R$ and
$F(\infty)=\lim_{y\to\infty}F(y)$.  
Hence, they
are left continuous on $(-\infty,\infty]$, vanish at $-\infty$ and 
can have a jump
discontinuity at $-\infty$ but not at $\infty$.
Then $H_{-\infty}\in\B$.  For $F\in\B$ and $g\in\bv$ define
$\langle F',g\rangle=\intinf g(x)\,dF(x)$.  Then $\langle H'_{-\infty},g\rangle
=g(-\infty)$.  Hence, for this type of differentiation with respect
to functions of bounded variation, $H_{-\infty}'\not=0$.  The Banach
lattice $\B$ then has $H_{-\infty}$ as an order unit.  We will explore
this type of differentiation elsewhere.

The usual pointwise ordering makes $L^1$ into a Banach lattice.
But the space of Henstock--Kurzweil integrable functions is not
a vector lattice. It is not closed under supremum and infimum since there
are functions integrable in this sense for which $\intinf f(x)\,dx$
converges but $\intinf |f(x)|\,dx$ diverges.  For example, 
$f(x)=x^2\sin(\exp(x^2))$.  Thus, even for functions,
when we allow conditional convergence we must look elsewhere to find
a lattice structure.

This
order $\preceq$ is not compatible with the usual order on distributions:
if $T,U\in\Dp$ then $T\geq U$ if and only if $\langle T-U,\phi\rangle
\geq 0$ for all $\phi\in\D$ such that $\phi\geq 0$.  If $T\geq 0$
then it is known that $T$ is a Borel measure.  The usual ordering on
distributions does
not give a vector lattice on $\acn$. 
For
example, if $F(x)=H_0(x)\int_0^x\sin(t)\,dt/t$ then $F(x)\geq 0$ 
for all $x\in\R$ so $F^{(n)}\succeq 0$ in
$\acn$.  With the distributional ordering, $\sup(F',0)$ is the function
equal to $\sin(x)/x$ when $x\in(2n\pi, (2n+1)\pi]$ for some integer
$n\geq 0$ and is equal to
$0$ otherwise.  This function is not in $\ac^1$ since the integral
defining $F$ converges conditionally.
None of the derivatives $F^{(n)}(x)$ are positive in
the pointwise  or
distributional sense.  Note that in $\acn$ we have
$(D^nF)^+=|F^{(n)}|=F^{(n)}$ and $(D^nF)^-=0$.

If two distributions are in more than one of the $\acn$ spaces
they may have different order relations in each such space.
For example, let $f(x)=-2(1-2x^2)\exp(-x^2)$, $F_1(x)=-2x\exp(-x^2)$
and $F_2(x)=\exp(-x^2)$.  Then $f, F_1, F_2\in \Bc$ and
we have the pointwise derivatives $f(x)=F_1'(x)=F_2''(x)$ for
each $x\in\R$.  Hence, $f\in\ac^1\cap\ac^2$.  Since $f(0)=-2<0$
and $f(1)=2/e>0$, it follows that $f$ is neither positive nor negative
in $\Bc$.
Since $F_1(-1)=2/e>0$ and $F_1(1)=-2/e<0$, it follows that $f$ is
neither positive nor negative in $\ac^1$. But $F_2(x)>0$ for all $x\in\R$
so $f\succeq 0$ in $\ac^2$.

A vector lattice is {\it order complete} (or {\it Dedekind complete})
if every nonempty subset that is bounded
above has a supremum.  But $\Br$ is not complete.
Let $F_n(x)=H_0(x-1/n)\sin(\pi/x)$ with $F_n(0)=0$.
Let $S=\{F_n\mid n\in\N\}$ then $S\subset\Bc$.
An upper bound for $S$ is the Heaviside step function $H_0$ 
but $\sup(S)(x)=H_0(x)\sin(\pi/x)$,
which is not regulated.
Hence, $\arn$ is also not complete.

A vector lattice is {\it Archimedean} if whenever $0\leq x\leq ny$ for
all $n\in\N$ and some $y\geq 0$ then $x=0$.  Applying the Archimedean
property at each point of $\R$ shows $\Br$ and hence $\arn$ are Archimedean.
All lattice inequalities that hold in $\R$ also hold in all Archimedean spaces
and all lattice equalities that hold in $\R$ also hold in all vector lattices.
See \cite{aliprantisborder}.  This expands the list of identities and
inequalities proved in Theorem~\ref{theoremlattice}.

A Banach lattice is an {\it abstract $L$ space} if
$\norm{x+y}=\norm{x}+\norm{y}$ for all $x,y\geq 0$.  
A Banach lattice is an {\it abstract $M$ space} if
$\norm{x\vee y}=\max(\norm{x},\norm{y})$ for all $x,y\geq 0$.  
See, for example \cite{aliprantisborder}.  We next
show that $\Br$ and $\arn$ are abstract $M$ spaces but neither
is an abstract $L$ space.
\begin{theorem}
(a) All of $\Br$, $\Bc$, $\arn$ and $\acn$ are abstract $M$ spaces.
(b) None of $\Br$, $\Bc$, $\arn$ or $\acn$ are abstract $L$ spaces.
\end{theorem}

\bigskip
\noindent
{\bf Proof:} It suffices to prove $\Br$ is an abstract $M$ space.  
(a) If $F_1, F_1\geq 0$ in $\Br$ then $\norm{F_1\vee F_2}_\infty
=\sup_{x\in\R}\max(F_1(x),F_2(x))\geq \sup_{x\in\R}F_1(x)=
\norm{F_1}_\infty$.  Similarly, $\norm{F_1\vee F_2}_\infty\geq
\norm{F_2}_\infty$.  So, $\norm{F_1\vee F_2}_\infty\geq
\max(\norm{F_1}_\infty,
\norm{F_2}_\infty)$.  And, $\norm{F_1\vee F_2}_\infty=
\sup_{x\in\R}\max(F_1(x), F_2(x))\leq \sup_{x\in\R}
\max(\norm{F_1}_\infty,\norm{F_2}_\infty)=\max(\norm{F_1}_\infty,
\norm{F_2}_\infty)$.    Hence, $\Br$ is an abstract $M$ space.

(b) It suffices to show $\Bc$ is not an abstract $L$ space.
Let $F_1(x)=1-|x|$ for $|x|\leq 1$ and $F_1(x)=0$, otherwise.
Let $F_2(x)=1-|x-2|$ for $|x-2|\leq 1$ and $F_2(x)=0$, otherwise.
Then $F_1(x), F_2(x)\geq 0$ for all $x\in\Rbar$.  And,
$\norm{F_1}_\infty=\norm{F_2}_\infty=\norm{F_1+F_2}_\infty=1$.  So
$\Bc$ is not an abstract $L$ space.\qed

For every measure $\mu$ it is
known that $L^1(\mu)$ is an abstract $L$ space and that a Banach
lattice is an abstract $L$ space if and only if it is lattice
isometric to  
$L^1(\nu)$ for some
measure $\nu$.  A Banach lattice is an abstract $M$ space with unit  
if and only if it is lattice isometric to $C(K)$ for some compact
Hausdorff space $K$.  The space $C(K)$ is the set of all real-valued
continuous functions on $K$.  These results are due to S.~Kakutani, M.~Krein and
others.  For references see \cite{aliprantisborder}.
The fact that $\acn$ is an abstract $M$ space but not an abstract
$L$ space indicates that what we have termed an integral here is
fundamentally different from the Lebesgue integral.  

\section{Banach algebra}\label{sectionalgebra}
A {\it commutative algebra} is a vector space $V$ over scalar field $\R$
with a multiplication $V\times V\mapsto V$ such that for all
$u,v,w\in V$ and all $a\in\R$, $u(vw) =(uv)w$ (associative),
$uv=vu$ (commutative),
$u(v+w)=uv+uw$ and $(u+v)w=uw+vw$ (distributive),
$a(uv)=(au)v$. 
If $(V,\norm{\cdot}_V)$ is a Banach space and $\norm{uv}_V\leq
\norm{u}_V\norm{v}_V$ then it is a Banach algebra.
For any compact Hausdorff space, $K$, the set of continuous
real-valued functions $C(K)$ is a commutative Banach algebra
under pointwise multiplication and
the uniform norm.  Since $\Rbar$ is compact and $\Br$ and $\Bc$ are
closed under pointwise multiplication, $\Br$ is a subalgebra of
$C(\Rbar)$ and $\Bc$ is a subalgebra of $\Br$.
The usual pointwise multiplication, $(FG)(x)=[F(x)][G(x)]$
for all $x\in\Rbar$, then makes $\Br$ into a commutative algebra.
The inequality
$\norm{F_1F_2}_\infty\leq\norm{F_1}_\infty\norm{F_2}_\infty$ for
all $F_1,F_2\in\Br$ shows
$\Br$ is a commutative Banach algebra.  

There is no unit.  For suppose $F(x)>0$ for all $x\in\R$.  If $eF=F$ 
then $e(x)=1$ for all $x\in\R$ so $e\not\in\Br$.
Consider the sequence $(u_n)\subset\Bc$ defined by $u_n(x)=
0$ for $x\leq -n$, $u_n(x)=x+n$ for $-n\leq x\leq 1-n$ and
$u_n(x)=1$ for $x\geq 1-n$.  For each $F\in\Bc$ we have
$\norm{F-u_nF}_\infty\to 0$.  Given $\epsilon>0$ there is $a\in\R$
such that $|F(x)|<\epsilon$ for all $x\leq a$.  We then have
$|F(x)-u_n(x)F(x)|=|F(x)||1-u_n(x)|<\epsilon$ for all $x\leq a$.
If $x\geq a$ take $n\geq 1-a$.  Then $u_n(x)=1$.  Hence,
$\norm{F-u_nF}_\infty\to 0$.  $\Bc$ is then said to have an
{\it approximate identity}.

By Lemma~\ref{lemma}, $\arn$ is a commutative Banach algebra, isomorphic
to $\Br$ for each $n\in\N$.  If $f_1, f_2\in\arn$ with respective
primitives $F_1,F_2\in\Br$ then $f_1f_2=D^n(F_1F_2)$.
\begin{theorem}
For each $n\in\N$, $\arn$ is a commutative Banach algebra without
unit, with approximate identity, isomorphic to $\Br$.
Similarly with $\acn$ and $\Bc$.
\end{theorem}

There is no difficulty in allowing functions in $\Br$ to be
complex-valued and using $\C$ as the field of scalars.  
Complex conjugation is then an involution on $\Br$.  Then
$\Br$ is a $C^\ast$-algebra since for each $F\in\Br$ we have
$\norm{{\overline F}}_\infty =\norm{F}_\infty$ and
$\norm{F{\overline F}}_\infty=\norm{F}_\infty^2$.  Thus, each
space $\arn$ is also a  $C^\ast$-algebra.

If $f_1,f_2\in\arn$ with respective primitives $F_1,F_2\in\Br$
then for all $h\in\Ibv^{n-1}$, 
\begin{align*}
&\intinf (f_1f_2)h  =  (-1)^{n-1}\intinf
D^n(F_1F_2) D^{n-1}h\\
&\qquad = (-1)^{n-1}F_1(\infty)F_2(\infty)h^{(n-1)}(\infty)
-(-1)^{n-1}\intinf F_1(x)F_2(x)\,dh^{(n-1)}(x).
\end{align*}
Let $a<\infty$ and 
$h=I^{n-1}[\chi_{(-\infty,a)}]$ with $\lambda=1$.
Then $\intinf (f_1f_2)h=(-1)^{n}F_1(a)F_2(a)$.
In particular, in $\ac^1$ we have $\int_{-\infty}^a(f_1f_2)=(\int_{-\infty}^a
f_1)(\int_{-\infty}^af_2)$.

There are zero divisors.  Let $F_1, F_2\in\D$ with disjoint
supports.  Then $F_1F_2=0$ in $\Br$ so $F_1^{(n)}F_2^{(n)}=0$ in
$\arn$, yet neither $F_1^{(n)}$ nor $F_2^{(n)}$ need be zero.
This example also shows the multiplication introduced in $\arn$
is not compatible with pointwise multiplication in the case when
elements of $\arn$ are functions.

If two distributions are in more than one of the $\acn$ spaces
they may have a different product in each such space.  For example,
let $f(x)=\sin(x)$ for $|x|\leq \pi$ and $f(x)=0$, otherwise.  
Then $f\in\Bc$ and its square in $\Bc$ is the function $f^2(x)=
\sin^2(x)$ for $|x|\leq \pi$ and $f^2(x)=0$, otherwise.  Now let
$F_1(x)=-1-\cos(x)$ for $|x|\leq \pi$ and $F_1(x)=0$, otherwise.
Let 
$F_2(x)=0$ for $x\leq-\pi$,
$F_2(x)=-x-\sin(x)-\pi$ for $|x|\leq \pi$ and $F_2(x)=-2\pi$ for
$x\geq \pi$.
Then $F_1,F_2\in\Bc$ and we have the pointwise derivatives
$f(x)=F_1'(x)=F_2''(x)$ for each $x\in\R$.  Hence, $f\in\ac^1\cap\ac^2$.
In $\ac^1$, $f^2(x)=D(F_1^2)(x)=-2[1+\cos(x)]\sin(x)$
for $|x|\leq \pi$ and $f^2(x)=0$, otherwise. 
In $\ac^2$, $f^2(x)=D^2(F_2^2)(x)=2[(1+\cos(x))^2-(x+\sin(x)+\pi)\sin(x)]$
for $|x|\leq \pi$ and $f^2(x)=0$, otherwise.  Hence, $f\in\Bc\cap\ac^1\cap
\ac^2$ but has a different product in each of these three spaces. 

By Proposition~\ref{propdelta}, for each integer $n\geq 0$, 
$\delta^{(n)}\in\ar^{n+1}$.  Let $\lambda\in\R$.  Then
$\delta^{(n)}\delta^{(n)}=D^{n+1}[H_\lambda^2]=D^{n+1}[H_{\lambda^2}]=
\delta^{(n)}$.

It is easy to see that $\Bc$ is a maximal ideal of $C(\Rbar)$.  See
\cite{kaniuth} for the definition.  It then follows that $\acn$ is a
maximal ideal of the space of $n$th derivatives of functions in $C(\Rbar)$.
Similarly, $\Br$ is a maximal ideal of the space $\B$, introduced in
Section~\ref{sectionlattice}, consisting of the left continuous
regulated functions
on $\Rbar$
such that $F(-\infty)=0$.  It then follows that $\arn$ is a maximal
ideal of the $n$th derivatives of such functions.  Note that
$\B$ has a unit, $H_{-\infty}=\chi_{(-\infty,\infty]}$.  As pointed
out in Section~\ref{sectionlattice}, the distributional derivative is
too coarse to distinguish between $H_{-\infty}$ and the constant
functions so this would entail using a finer notion of derivative
such as using functions of bounded variation for test functions.
This is something we will pursue elsewhere.

Define $AC_{-\infty}(\Rbar)$ to be the
functions in $AC(\Rbar)$ whose limit vanishes at $-\infty$.  
And, $AC_{-\infty}^n(\R)$ consists of the distributions $f\in\Dp$ such that
$f=D^nF$ for some $F\in AC_{-\infty}(\Rbar)$.   Note that $AC_{-\infty}^1(\Rbar)
=L^1$.
It is easy to show that $AC_{-\infty}(\Rbar)$ is closed under pointwise
multiplication.  Hence, it is a Banach subalgebra of $\Bc$.  Then
$AC_{-\infty}^n(\R)$ is a Banach subalgebra of $\acn$.  Of course,
$L^1$ is not an algebra under the usual pointwise multiplication.  
Similarly for the
spaces of $n$th derivatives of primitives of Henstock--Kurzweil and wide Denjoy integrable
functions.  See \cite{gordon} for the definitions of the relevant spaces
of primitives.  Under convolution $L^1$ is a Banach algebra.  Although
convolution has been defined in $\ac^1\times L^1$ in \cite{talvilaconvolution}
it does not seem possible to define convolution in $\ac^1\times\ac^1$.
Convolutions can be defined for distributions but restrictions on the
supports are generally imposed.  See \cite{zemanian}.

\end{document}